\documentclass[a4,12pt]{amsart}%

\usepackage{amsfonts,amssymb,amsmath,amsthm}
\usepackage{algorithmicx, algorithm}
\usepackage{algcompatible}
\usepackage{color,ulem}
\usepackage{fullpage}
\usepackage{graphicx}%
\usepackage[active]{srcltx}
\usepackage{mathabx}
\usepackage{verbatim}
\usepackage{cancel}
\usepackage{mathtools}
\usepackage{tikz,pgfplots}
\usetikzlibrary{snakes}
\usepackage{cite}
\usepackage{hyperref}
\hypersetup{
  linkcolor  = green!60!black,
  citecolor  = blue!60!black,
  urlcolor   = red!60!black,
  colorlinks = true,
}

\usepackage{caption}
\usepackage{subcaption}

\usepackage[official]{eurosym}

\providecommand{\U}[1]{\protect\rule{.1in}{.1in}}

\newtheorem{theorem}{Theorem}

\newtheorem{definition}[theorem]{Definition}

\newtheorem{proposition}[theorem]{Proposition}
\newtheorem{remark}[theorem]{Remark}

\def\d{\hbox{d}}

\newcommand{\lm}[1]{{\textcolor{black}{#1}}}
\newcommand{\lmxxx}[1]{{\textcolor{black}{#1}}}
\usepackage{lipsum}
\makeatletter
\g@addto@macro{\endabstract}{\@setabstract}
\newcommand{\authorfootnotes}{\renewcommand\thefootnote{\@fnsymbol\c@footnote}}%
\makeatother

\begin{document}
\title{A \lmxxx{Feynman-Kac formula} approach for computing \lmxxx{expectations and threshold crossing probabilities} of non-smooth stochastic dynamical systems}
\maketitle
\begin{center}
  \normalsize
  \authorfootnotes
  Laurent Mertz\footnote{ECNU-NYU Institute of Mathematical Sciences,
    NYU Shanghai, Shanghai, China},
  Georg Stadler\footnote{Courant Institute of Mathematical Sciences,
    New York University, New York, USA},
  \and
  Jonathan Wylie\footnote{Department of Mathematics, City University
    of Hong Kong, Hong Kong, China}
\end{center}

\begin{abstract}
  We present a computational alternative to probabilistic simulations
  for non-smooth stochastic dynamical systems that are prevalent
  in engineering mechanics. As examples, we target (1) stochastic
  elasto-plastic problems, which involve transitions between elastic and
  plastic states, and (2) obstacle problems with noise, which involve
  discrete impulses due to collisions with an obstacle. We formally introduce a class of partial differential equations related to the Feynman-Kac formula, where the underlying stochastic processes satisfy variational inequalities modelling elasto-plastic and obstacle oscillators. 
  \lmxxx{We then focus on solving them numerically.} The main challenge in solving \lmxxx{these equations} is the non-standard boundary conditions
  which describe the behavior of the underlying process on the
  boundary. \lmxxx{We illustrate how to use our approach to compute expectations and other statistical quantities, such as the asymptotic growth rate of variance in asymptotic formulae for threshold crossing probabilities.}
\end{abstract}

\vspace{5pt}

\section{Motivations and Goal}
\noindent In this paper, we propose a computational alternative to probabilistic simulations for a certain type of non-smooth stochastic dynamical systems, namely
\begin{equation}
\label{Nonsmooth}
\dot{X}_t = F(X_t,Y_t),
\quad \dot{Y}_t = G(X_t,Y_t) + \lm{\mathbf{H}_t} + \dot{W}_t, \quad \forall t > 0
\end{equation}
with initial condition $(X_0,Y_0) = (x,y) \in \mathbb{R}^2$. Here, $(X_t, Y_t) \in \mathbb{R}^2$ is the state variable at time $t$, dots denote derivatives with respect to time, $F(x,y)$ and $G(x,y)$ are deterministic functions, $\dot{W}$ represents a white noise random forcing in the sense that $W$ is a Wiener process, and $\mathbf{H}$ is a functional depending \lm{on interactions of the state variable at boundaries and interfaces}. We are interested in statistical quantities that characterize and predict the behavior of $\{ (X_t,Y_t), t \geq 0 \}$ such as means, moments, correlations \lmxxx{and probabilities of crossing thresholds}. These are of the form
\[
A 
\triangleq 
\mathbb{E} \left ( f(X_T,Y_T) + \int_0^T g(X_\tau,Y_\tau) \d \tau,  \right ),
\quad
B 
\triangleq 
\mathbb{E} \left ( \int_0^\infty e^{-\lambda \tau} g(X_\tau,Y_\tau) \d \tau \right ),
\]
\[
C \triangleq \lim_{T \to \infty}  \mathbb{E} \left ( f(X_T,Y_T) \right ), 
\quad
\lmxxx{D \triangleq \mathbb{P} \left ( \max_{0 \leq t \leq T} \left | \int_0^t f(X_\tau,Y_\tau) \d \tau \right | \geq b \right ),} 
\]
and
\begin{align*}
A' & \triangleq \mathbb{E} \left [ \left ( f(X_T,Y_T) + \int_0^T g(X_\tau,Y_\tau) \d \tau  \right ) \left ( \varphi(X_{T+h},Y_{T+h}) + \int_0^{T+h} \psi(X_\tau,Y_\tau) \d \tau,  \right )  \right ],\\
B' & \triangleq \mathbb{E} \left (  \int_0^\infty  \int_0^\infty e^{-\lambda \tau - \mu \theta} g(X_\tau,Y_\tau) \psi (X_\theta,Y_\theta) \d \tau \d \theta \right ),\\ 
C' & \triangleq \lim_{T \to \infty} \frac{1}{T} \mathbb{E} \left ( \int_0^T \int_0^T g(X_\tau,Y_\tau) \psi (X_\theta,Y_\theta) \d \tau \d \theta \right ). 
\end{align*}
Here and in the remainder of the paper, $\mu, \lambda$ are positive
numbers, \lmxxx{$b$ is a given threshold}, $T>0$ is a given time, $h \geq 0$ and $f,g, \varphi, \psi$
are continuous functions. For cases in which $\mathbf{H} \equiv 0$ and
$(F,G)$ satisfy appropriate smoothness conditions, a natural setting
for characterizing such quantities \lmxxx{with partial differential equations (PDEs) is to use the Feynman-Kac formula (FKf) \cite{R96}.} 
Furthermore, if in addition $F$ and $G$
can be written in terms of a Hamiltonian structure $\mathcal{H}(x,y)$,
in the sense
\[
F(x,y) \triangleq \frac{\partial \mathcal{H}}{\partial x}(x,y), \quad  G(x,y) \triangleq \frac{\partial \mathcal{H}}{\partial y}(x,y) + \eta(x,y) \frac{\partial \mathcal{H}}{\partial x}(x,y), 
\quad (x,y) \in \mathbb{R}^2,
\] 
where $\eta$ is a well behaved function, then the process $(X_t,Y_t)$ belongs to the class of Stochastic Hamiltonian Systems (SHS) \cite{T02} for which quantities of type $C$ and $C'$ are well defined. Applications that make use of this framework abound in many areas of science and engineering, e.g., finance, chemistry, biology, neuroscience, economy and mechanics.\\

\noindent There are many important applications
involving non-smooth dynamics in which quantities 
of the form $A,B,C, \lmxxx{D}$ and $A',B',C'$ are of interest. 
We have in mind problems involving interactions with boundaries,
constraints, phase transitions or hysteresis. An important class of
examples includes elasto-plasticity (EP) problems with random forcing
\cite{KS66,F08} where the dynamics takes into account phase
transitions between elastic and plastic states. A second class of
examples includes stochastically driven obstacle problems
\cite{B98,DI04} where the dynamics has to take into account instantaneous collisions
with an obstacle. The constitutive models are shown in Figure~\ref{fig1}.
A number of authors have applied innovative techniques to determine statistics of non-smooth oscillators
\cite{TD01, DB93, YBDFP17, XWJHY13, Z14, RWXF10, RS03,NM96,KYS94, KS09}.
They used methods such as stochastic averaging, statistical linearization and numerical path integration.
However, some of these methods are based on approximations of the underlying mechanical system, whereas the methodology 
that we present in this paper does not make any such approximations and therefore provides a rigorous basis for calculating
the above-mentioned quantities. A concise review on the stochastic engineering dynamics literature is provided in section 2.\\

{
\color{black}
\begin{figure}[htbp]
  \centering
  \begin{tikzpicture}[scale=1.0]
\draw[thick] (-1.0,1.00) node[right] {$\bold{(a)}$};
\draw[color=gray] (-1.125,0.875) -- (-1,1);
\draw[color=gray] (-1.125,0.75) -- (-1,0.875);
\draw[color=gray] (-1.125,0.625) -- (-1,0.75);
\draw[color=gray] (-1.125,0.5) -- (-1,0.625);
\draw[color=gray] (-1.125,0.375) -- (-1,0.5);
\draw[color=gray] (-1.125,0.25) -- (-1,0.375);
\draw[color=gray] (-1.125,0.125) -- (-1,0.25);
\draw[color=gray] (-1.125,0) -- (-1,0.125);
\draw[color=gray] (-1.125,0.125) -- (-1,0.25);
\draw[color=gray] (-1.125,-0.125) -- (-1,0);
\draw[color=gray] (-1.125,-0.25) -- (-1,-0.125);
\draw[color=gray] (-1.125,-0.375) -- (-1,-0.25);
\draw[color=gray] (-1.125,-0.5) -- (-1,-0.375);
\draw[color=gray] (-1.125,-0.625) -- (-1,-0.5);
\draw[color=gray] (-1.125,-0.75) -- (-1,-0.625);
\draw[color=gray] (-1.125,-0.875) -- (-1,-0.75);
\draw[color=gray] (-1.125,-1) -- (-1,-0.875);
\draw[thick] (-1,-1) -- (-1,1) ;
\draw[thick] (-1,0.5) -- (0.125,0.5) ;
\draw[thick] (0.25,0.5) -- (1,0.5) ; 
\draw[thick] (-0.25,0.625) -- (0.25,0.625) ;
\draw[thick] (-0.25,0.375) -- (0.25,0.375) ;
\draw[thick] (0.25,0.375) -- (0.25,0.625) ;
\draw[thick] (0.125,0.4) -- (0.125,0.6) ;
\draw[thick] (-1,-.5) -- (-0.9,-.5) ;
\draw[snake=coil,segment length=4pt, thick] (-0.90,-.5) -- (0.1,-.5) ; 
\draw[thick] (0.1,-0.5) -- (0.5,-0.5) ;
\draw[thick] (0.4,-0.375) -- (0.9,-0.375) ;
\draw[thick,->] (0.4,-0.375) -- (0.4,-0.5);
\draw[thick] (0.9,-0.375) -- (0.9,-0.5);
\draw[thick] (.8,-0.5) -- (1,-0.5) ;
\draw[step=0.0125cm,color=black] (1.5,-0.125) grid (1.75,0.125);
\draw[thick] (1,-0.5) -- (1,0.5);
\draw[thick] (1,0) -- (1.5,0) ;
\draw[->]  (2.25,0) -- (1.80,0);
\draw[thick] (2.25,0) node[right] {$\dot W_t$};
\end{tikzpicture}
\hspace{1cm}
  \begin{tikzpicture}[scale=1.0]
\draw[thick] (-1.0,1.00) node[right] {$\bold{(b)}$};
\draw[color=gray] (-1.125,0.875) -- (-1,1);
\draw[color=gray] (-1.125,0.75) -- (-1,0.875);
\draw[color=gray] (-1.125,0.625) -- (-1,0.75);
\draw[color=gray] (-1.125,0.5) -- (-1,0.625);
\draw[color=gray] (-1.125,0.375) -- (-1,0.5);
\draw[color=gray] (-1.125,0.25) -- (-1,0.375);
\draw[color=gray] (-1.125,0.125) -- (-1,0.25);
\draw[color=gray] (-1.125,0) -- (-1,0.125);
\draw[color=gray] (-1.125,0.125) -- (-1,0.25);
\draw[color=gray] (-1.125,-0.125) -- (-1,0);
\draw[color=gray] (-1.125,-0.25) -- (-1,-0.125);
\draw[color=gray] (-1.125,-0.375) -- (-1,-0.25);
\draw[color=gray] (-1.125,-0.5) -- (-1,-0.375);
\draw[color=gray] (-1.125,-0.625) -- (-1,-0.5);
\draw[color=gray] (-1.125,-0.75) -- (-1,-0.625);
\draw[color=gray] (-1.125,-0.875) -- (-1,-0.75);
\draw[color=gray] (-1.125,-1) -- (-1,-0.875);
\draw[thick] (-1,-1) -- (-1,1) ;
\draw[thick] (-1,0.5) -- (0.125,0.5) ;
\draw[thick] (0.25,0.5) -- (1,0.5) ; 
\draw[thick] (-0.25,0.625) -- (0.25,0.625) ;
\draw[thick] (-0.25,0.375) -- (0.25,0.375) ;
\draw[thick] (0.25,0.375) -- (0.25,0.625) ;
\draw[thick] (0.125,0.4) -- (0.125,0.6) ;
\draw[thick] (-1,-0.5) -- (-0.5,-0.5) ;
\draw[snake=coil,segment length=4pt, thick] (-0.5,-0.5) -- (0.5,-0.5) ; 
\draw[thick] (0.5,-0.5) -- (1,-0.5) ;
\draw[step=0.0125cm,color=black] (1.5,-0.125) grid (1.75,0.125);
\draw[thick] (1,-0.5) -- (1,0.5);
\draw[thick] (1,0) -- (1.5,0) ;
\draw[->]  (2.25,0) -- (1.80,0);
\draw[thick] (2.25,0) node[right] {$\dot W_t$};

\draw[color=gray] (1.125,0.875) -- (1.25,1);
\draw[color=gray] (1.125,0.75) -- (1.25,0.875);
\draw[color=gray] (1.125,0.625) -- (1.25,0.75);
\draw[color=gray] (1.125,0.5) -- (1.25,0.625);
\draw[color=gray] (1.125,0.375) -- (1.25,0.5);
\draw[color=gray] (1.125,0.25) -- (1.25,0.375);
\draw[color=gray] (1.125,0.125) -- (1.25,0.25);
\draw[color=gray] (1.125,0) -- (1.25,0.125);
\draw[color=gray] (1.125,0.125) -- (1.25,0.25);
\draw[color=gray] (1.125,-0.125) -- (1.25,0);
\draw[color=gray] (1.125,-0.25) -- (1.25,-0.125);
\draw[color=gray] (1.125,-0.375) -- (1.25,-0.25);
\draw[color=gray] (1.125,-0.5) -- (1.25,-0.375);
\draw[color=gray] (1.125,-0.625) -- (1.25,-0.5);
\draw[color=gray] (1.125,-0.75) -- (1.25,-0.625);
\draw[color=gray] (1.125,-0.875) -- (1.25,-0.75);
\draw[color=gray] (1.125,-1) -- (1.25,-0.875);
\draw[thick] (1.25,-1) -- (1.25,1) ;
\draw[color=gray] (2,0.875) -- (2.125,1);
\draw[color=gray] (2,0.75) -- (2.125,0.875);
\draw[color=gray] (2,0.625) -- (2.125,0.75);
\draw[color=gray] (2,0.5) -- (2.125,0.625);
\draw[color=gray] (2,0.375) -- (2.125,0.5);
\draw[color=gray] (2,0.25) -- (2.125,0.375);
\draw[color=gray] (2,0.125) -- (2.125,0.25);
\draw[color=gray] (2,0) -- (2.125,0.125);
\draw[color=gray] (2,0.125) -- (2.125,0.25);
\draw[color=gray] (2,-0.125) -- (2.125,0);
\draw[color=gray] (2,-0.25) -- (2.125,-0.125);
\draw[color=gray] (2,-0.375) -- (2.125,-0.25);
\draw[color=gray] (2,-0.5) -- (2.125,-0.375);
\draw[color=gray] (2,-0.625) -- (2.125,-0.5);
\draw[color=gray] (2,-0.75) -- (2.125,-0.625);
\draw[color=gray] (2,-0.875) -- (2.125,-0.75);
\draw[color=gray] (2,-1) -- (2.125,-0.875);
\draw[thick] (2,-1) -- (2,1) ;
\end{tikzpicture}
\caption{\lmxxx{Constitutive models: \textbf{(a)} an elasto-perfectly-plastic oscillator and \textbf{(b)} an oscillator with an obstacle and impacts. A mass
(black box) is associated in series with elements which are themselves an association in parallel or in series of elementary rheological models. Each whole system is excited by a time-dependent random forcing $\dot W_t$.}}
\label{fig1}
\end{figure}
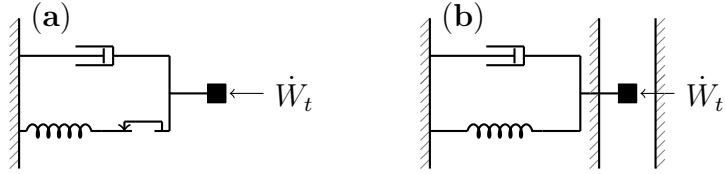
}
\noindent \lm{From a mathematical viewpoint, these problems can be reformulated in terms of (degenerate) stochastic variational inequalities (SVIs) \cite{BT08,BD06}}. With this formulation, the evolution of the state variables is Markovian and Kolmogorov Equations (KEs) can be derived.\\

\noindent 
Except for very few non-smooth dynamical systems, where analytical
expressions are available, in general one has to resort to
computational techniques to determine $A,B,C, \lmxxx{D}$ and $A',B',C'$. 
The most straightforward computational methods are direct probabilistic simulations. Such techniques are simple to
implement and widely used. For low-dimensional state variables, they
are less efficient than techniques based on partial differential equations (PDEs), 
provided the latter are available. For the two non-smooth stochastic processes targeted in this
paper, the presence of constraints or phase transitions (governed by a variational inequality structure) leads to
non-standard boundary conditions in the \lmxxx{PDEs for FKf}. These boundary conditions characterize the
behavior of the underlying process on the boundary. This paper is concerned with the
numerical treatment of such non-standard \lmxxx{PDEs} and their application to compute the quantities of $A,B,C\lmxxx{,D}$ and $A',B',C'$.\\

\noindent We point out that the idea of solving non-standard \lmxxx{PDEs for FKf} in the context of a SVI has been used before \cite{BMPT09}. However, the technique in \cite{BMPT09} is specific to the elasto-perfectly-plastic oscillator excited by white noise, for which the non-standard condition has to be satisfied in a finite (namely two) number of points. It has no natural extension to problems of obstacle type because for these problems the non-standard condition must be satisfied on a continuous set of points.\\

\noindent In the remainder of this section, we give the definitions of the 
elasto-plastic and obstacle models in presence of noise, and we discuss the quantities we are interested in.

\subsection{Elasto-plastic problem with noise} \label{epo} A basic
prototype for modeling mechanical structures that admit permanent
deformation under vibrations is the elastic-plastic oscillator
\cite{KS66}. The dynamics focuses on
two quantities: the total deformation, $X_t$, supported by the structure
when subjected to vibrations, and its velocity, $\dot{X}_t$. 
For general elasto-plastic problems, the nonlinear functional $\mathbf{H}$ in \eqref{Nonsmooth} is a restoring force arising from the structure. The exact form of $\mathbf{H}$ depends on the particular structure in question. 
The nonlinearity in such models comes from the switching of regimes from an elastic phase to a plastic one, or vice versa. 
For the elasto-perfectly-plastic oscillator (EPPO) model, the irreversible (plastic) deformation $\Delta$ and the reversible (elastic) deformation $Z$ at time $t$ satisfy  
\begin{eqnarray*}
\dot{Z}_t & = \dot{X}_t, \quad \dot{\Delta}_t & = 0, \quad \mbox{in elastic phase},\\
\dot{\Delta}_t & = \dot{X}_t, \quad \dot{Z}_t & = 0, \quad \mbox{in plastic phase},
\end{eqnarray*}
where $X_t = Z_t + \Delta_t$. Typically, $|Z_t|$ is bounded by a given threshold $P_Y$ at all times $t$, the system is in the plastic phase when $ \vert Z_t \vert = P_Y$ and the elastic phase when $\vert Z_t \vert < P_Y$. Here, $P_Y$ is an elasto-plastic bound, known as the ``Plastic Yield'' in the engineering literature. We assume the  force $\mathbf{H}$ is a linear function of $Z$ (while $\Delta$ remains constant) in the elastic phase and a constant (while $Z$ remains constant at $\pm P_Y$) in the plastic phase as follows:
\begin{equation}
\label{pathdependent_ep}
\lm{\mathbf{H}_t} = k Z_t , \quad  k > 0.
\end{equation}
The permanent (plastic) deformation at time $t$ can then be
written as
\[
\Delta_t = \int_0^t  \mathbf{1}_{\{ |Z_s| = P_Y \} } Y_s \hbox{d} s.
\] 
We consider the case of linearly damped spring for which $F(x,y) = y$ and $G(x,y) = -c_0 y$.

\subsubsection*{SVI framework}
\noindent It has been shown that the dynamics of such a nonlinear
oscillator can be described mathematically by means of SVIs \cite{BT08}.
The dynamics is then described by the pair $(Y_t,Z_t)$ that
satisfies
\begin{equation}
\label{svi_ep}
\forall t, \forall | \phi | \leq P_Y, \: \d Y_t  = - (c_0 Y_t + kZ_t) \d t + \d W_t, \quad (\d Z_t-Y_t \d t)(\phi - Z_t) \geq 0
\end{equation}
and appropriate initial conditions for $Y_0$ and $Z_0$ must be prescribed. \textcolor{black}{Here, whereas $Y_t = \dot{X}_t$, $X_t$ is not involved in the dynamics}.\\

\noindent Characterizing the statistics of the state
variable using \lmxxx{FKf} is interesting for engineering purposes. 
To illustrate the robustness and efficiency of our approach we will compute \lmxxx{five} important quantities in this paper:

\begin{itemize}
\item the probability of finding the system in the plastic state
\begin{equation}
\label{epp_plasticstate}
\tag{$E_1$}
\mathbb{P}\left(\left|Z_T \right|=P_Y\right) \: \mbox{for} \: T \geq 0, 
\quad \mbox{and} \quad 
\lim_{T \to \infty} \mathbb{P}\left(\left|Z_T \right|=P_Y\right),
\end{equation}
\item the mean kinetic energy
\begin{equation}
\label{epp_mkenergy}
\tag{$E_2$}
\mathbb{E}Y_T^2 \: \mbox{for} \: T \geq 0,  
\quad \mbox{and} \quad \lim_{T \to \infty} \mathbb{E}Y_T^2,
\end{equation}
\item the variance of both the plastic deformation and the total deformation 
\begin{equation}
\label{epp_variance}
\tag{$E_3$}
\sigma^2\left(\Delta_T \right), \quad \sigma^2\left(X_T \right)\: \mbox{for} \: T \geq 0, 
\quad \mbox{and} \quad 
\lim_{T \to \infty}
\frac{1}{T}\sigma^2\left(\Delta_T \right), \quad \lim_{T \to \infty}
\frac{1}{T}\sigma^2\left(X_T \right), 
\end{equation}
\item
any correlation structure between $(Y_T,Z_T)$ and $(Y_{T+h},Z_{T+h})$ of the form 
\[
\mathbb{E} f(Z_T,Y_T) \varphi(Z_{T+h},Y_{T+h}),
\] 
where here we will target
\begin{equation}
\label{epp_correlation}
\tag{$E_4$}
\mathbb{P} \left ( | Z_{T+h}| = P_Y, |Z_T| = P_Y \right ), \quad \mathbb{E} Y_T^2 Y_{T+h}^2 \: \mbox{for} \: T \geq 0, 
\end{equation}
and
\[
\lim_{T \to \infty} \mathbb{P} \left ( | Z_{T+h}| = P_Y, |Z_T| = P_Y \right ), \quad \lim_{T \to \infty} \mathbb{E} Y_T^2 Y_{T+h}^2,
\]
\item
\lmxxx{
probabilities of threshold crossings for the total deformation
\begin{equation}
\label{epp_crossings}
\tag{$E_5$}
\mathbb{P} \left ( \max_{0 \leq t \leq T} | X_t | \geq b \right ), \: b, T \: \textup{large}.
\end{equation}
}
\end{itemize}
The probability of finding the system in the plastic state and the mean kinetic energy, at a given time $T$, are quantities of type $A$ whereas, at large time, they are of type $C$. Also, their Laplace transforms belong to type $B$. The variance of the plastic deformation and of the total deformation, at a given time $T$, are quantities of type $A'$ and their asymptotic growth belongs to the type $C'$. Any correlation structure of $(Y,Z)$ between two different instants $T$ and $T+h$ is of type $A'$ and its corresponding double Laplace transform (with parameters $\lambda$ and $\mu$) is of type $B'$.
\lmxxx{The probabilities of threshold crossings for the plastic deformation are of type $D$.}
 These quantities are challenging to compute because analytic formulas are not available and probabilistic simulations are computationally expensive.
 
\subsection{Obstacle problem with noise} 
\label{vio}
\noindent
One of the simplest models exhibiting a vibro-impact motion can be
expressed as a single degree of freedom oscillator constrained by
obstacles \cite{B98} located at $|X|=P_O$ (position of the obstacle).  It is common in the engineering literature to
formulate the dynamics of a stochastic obstacle oscillator in terms of
a stochastic process $X_t$, the oscillator displacement. For general obstacle problems, the nonlinearity in such models comes from the collisions.
If at a time $t$, the state hits the obstacle with incoming velocity $\dot{X}_t$, it immediately bounces
back with velocity $-e\dot{X}_t$, that is, $\dot{X}_{t+} = - e \dot{X}_{t-}$, $e \in [0,1]$.
The nonlinear functional $\mathbf{H}$ in \eqref{Nonsmooth} is a force keeping track of the past discrete impulses due to collisions with an obstacle, here
\begin{equation}
\label{pathdependent_op}
\lm{\mathbf{H}_t = \dot{I}_t.}
\end{equation}
The net impulse process, which keeps track of the sum of all past impulses, at time $t$ can be
written as
\[
I_t \triangleq \sum_{0 \leq s \leq t} \left ( \dot{X}_{s+} - \dot{X}_{s-} \right ) \mathbf{1}_{ \{ |X_s|=P_O \}}.
\] 
For simplicity, we consider the case of a linearly damped spring: $F(x,y) = y$ and $G(x,y) = - c_0 y - k x$.

\subsubsection*{Reflected Langevin process framework}
{
\color{black}
\noindent From a mathematical viewpoint (using a stochastic differential framework), the dynamics of such a nonlinear
oscillator can be described in the framework of a Reflected Langevin process \cite{J12,J13} of the form
\begin{equation}
X_t = x + \int_0^t Y_s \d s,
\quad
Y_t = y - \int_0^t ( c_0 Y_s + k X_s )\d s + W_t - (1+e) \sum_{0 \leq s \leq t} \dot{Y}_{s-} \mathbf{1}_{ \{ |X_s|=P_O \}} 
\label{svi_obstacle}
\end{equation}
where $Y_t \triangleq \dot{X}_t$.
\begin{remark}
\textcolor{black}{When $e=1$,
Equation \eqref{svi_obstacle} 
can be formulated in the framework of SVIs \cite{BD06} as follows
$
\forall t, \: \vert X_t \vert \leq  P_O, \: \forall  \vert \varphi \vert   \leq  P_O, \:  
(\d Y_t + (c_0 Y_t + k X_t) \d t - \d W_t)(\varphi - X_t)  \geq   0.
$}
\end{remark}
}
{
\color{black}
\subsubsection*{Comments}
The implementation of the impact rule can be explained using stopping times. We start by defining $\tau_0 \triangleq 0$ and $\{ (X_t^0,Y_t^0), t \geq 0 \}$ to be the solution of the unconstrained problem
\begin{equation}
\label{shs}
d X^0 = Y^0 dt, 
\quad 
d Y^0 = -(c_0 Y^0+k X^0) d t + d W 
\end{equation}
and define $\tau_1 \triangleq \inf \{ t > \tau_0, \: |X_t^0| = P_O \}$. For $t \geq \tau_1$, we define $\{ (X_t^1,Y_t^1), t \geq \tau_1 \}$ to be the solution of \eqref{shs}
with the initial condition $X_{\tau_1}^1 \triangleq X_{\tau_1}^0$ and $Y_{\tau_1}^1 \triangleq - e Y_{\tau_1}^0$ and then we define $\tau_2 \triangleq \inf \{ t > \tau_1, \: |X_t^1| = P_O \}$.  
Knowing $\tau_n,X^n,Y^n$, we can recursively define $\tau_{n+1},X^{n+1},Y^{n+1}$.
Thus, the whole process can be defined on each interval $[\tau_n,\tau_{n+1})$ by setting $(X,Y) \triangleq (X^n,Y^n)$.
It should be pointed out that the presentation remains formal as there are mathematical subtleties behind the definition of these stopping times.
An existence and uniqueness result in a weak sense (see for instance, definition 3.1 page 300 in \cite{KS91}) for an equation similar to \eqref{svi_obstacle} can be found in the works of Jacob \cite{J12,J13}. To be more precise, it concerns an integrated Wiener process constrained to stay in $[0,\infty)$ by a partially elastic boundary at $0$.
As explained in Jacob's works, even if the aforementioned process enjoys the property of local pathwise existence and uniqueness away from the obstacle, global pathwise uniqueness results do not hold. Indeed, even for the deterministic problem it has been shown in the work of Ballard \cite{B00} that the only case in which there exists a unique solution is when the forcing is an analytic function (locally given by a convergent power series), and this condition cannot be relaxed.
}\\

\noindent For the obstacle problem, we consider the following \lmxxx{five} quantities:
\begin{itemize}
\item
the probability of finding the system in the neighborhood of the obstacle with a low velocity, 
\begin{equation}
\label{op_proba}
\tag{$E_1'$}
\mathbb{P} \left ( f(X_T,Y_T) \leq \epsilon \right ), \: \mbox{for} \: T \geq 0, 
\quad \mbox{ and } \quad
\lim_{T \to \infty} \mathbb{P} \left ( f(X_T,Y_T) \leq \epsilon  \right )
\mbox{ for some } \epsilon >0
\end{equation}
where $f(x,y) \triangleq \sqrt{(|x| - P_O)^2 + y^2}$.
\item the mean kinetic energy
\begin{equation}
\label{op_mkenergy}
\tag{$E_2'$}
\mathbb{E}Y_T^2 \: \mbox{for} \: T \geq 0,  \quad \mbox{and} \quad \lim_{T \to \infty} \mathbb{E}Y_T^2,
\end{equation}
\item the variance of the integral of the displacement
\begin{equation}
\label{op_variance}
\tag{$E_3'$}
\sigma^2\left(\int_0^T X_s \d s\right) \: \mbox{for} \: T \geq 0,  
\quad \mbox{and} \quad \lim_{T \to \infty} \frac{1}{T} \sigma^2\left(\int_0^T X_s \d s\right).
\end{equation}
Up to a multiplicative constant, this is equivalent to computing the variance of the change in momentum due to the restoring force $kX$.
\item
any correlation structure between $(X_T,Y_T)$ and $(X_{T+h},Y_{T+h})$ of the form 
\begin{equation}
\label{op_correlation}
\tag{$E_4'$}
\mathbb{P} \left ( f(X_{T+h},Y_{T+h}) \leq \epsilon , f(X_T,Y_T) \leq \epsilon \right ), 
\quad \mathbb{E} Y_T^2 Y_{T+h}^2 \: \mbox{for} \: T \geq 0, 
\end{equation}
and
\[
\label{CorrStructObst2}
\lim_{T \to \infty} \mathbb{P} \left ( f(X_{T+h},Y_{T+h}) \leq \epsilon , f(X_T,Y_T) \leq \epsilon \right ),
\quad \lim_{T \to \infty} \mathbb{E} Y_T^2 Y_{T+h}^2,
\]
\item
\lmxxx{
probabilities of threshold crossings for the integral of the displacement
\begin{equation}
\label{op_crossings}
\tag{$E_5'$}
\mathbb{P} \left ( \max_{0 \leq t \leq T} \left | \int_0^t X_s d s \right | \geq b \right ), \: b, T \: \textup{large}.
\end{equation}
}
\end{itemize}
An explanation similar to what was provided above for the elasto-plastic problem (in terms of plastic state, variance of the plastic deformation and correlations) applies to the obstacle problem (in terms of mean kinetic energy, variance of the integral of the displacement, correlations \lmxxx{and probabilities of threshold crossings for the integral of the displacement}).

\subsection{Approach and overview} 
The PDEs related to $A,B,C$ and $A',B',C'$ including quantities $(E_i)$, $(E_i')$, $1 \leq i \leq 4$
for the elasto-plastic and obstacle problems are non-standard boundary value problems. 
For this type of non-standard PDEs, (for instance see Section 3.2), only partial
existence and uniqueness results are available, mainly for the case of
the EPPO problem with noise \cite{BT08,BT10, BMY16,BM15}. This is
because standard PDE theory techniques do not apply due to the
non-standard boundary conditions and the degeneracy of these
problems. Therefore, in this essay, we mostly study the behavior of
these PDEs numerically in order to gather insight in the solution
behavior and we compare with probabilistic approximations in order to
conjecture whether solutions exist or not. For this purpose, we solve
the PDE problems on sequences of grids with increasing resolution and
monitor the behavior of the numerical solutions. If the numerical
solutions computed from differently accurate discretization yield a
converging behavior as the mesh is refined, we can conjecture that the
continuous problem has a solution.\\

\noindent Monte Carlo methods are used to compare the PDE results with
probabilistic results. These comparisons increase our confidence in
the solution of the PDE problem, allow us to study approximation
errors in both schemes, and establish a direct connection between
the SVIs and the PDE problems.

\subsection{Organization of the paper} 
\textcolor{black}{A review of established stochastic engineering methodologies is given in section 2.}
In section 3, we show the connection between the quantities $A$ and $A'$ to the solution of non-standard \textit{parabolic problems},
the functions $B$ and $B'$ to the solution of non-standard \textit{elliptic problems} and the functions $C$ and $C'$ to the solution of non-standard \textit{Poisson problems}.
\lmxxx{In Section 4, an asymptotic formula for D is presented when $T$ and $b$ are large enough. This formula relies on quantities of type $C'$.}
In Section 5, we present a numerical approach for solving non-standard \lmxxx{PDEs}. The method is first presented and applied to the
elasto-plastic problem. Then, it is applied to the obstacle problem. Numerical results are presented.  
In Section 6, we compare the approach proposed in this paper with  previous techniques employed for a white noise EPPO.
Finally, in Section 7, broader impacts of the present method are discussed for promising engineering applications.

\section{Review of stochastic engineering methodologies}
{ 
\color{black}
\noindent We briefly review three main established techniques, namely i) stochastic averaging, ii) statistical linearization and iii) numerical path integration.
The notation used in this section only applies to this section and not beyond.

\subsection{Statistical linearization}
The statistical linearization (SL) method, proposed by \cite{C63}, has proved to be a very useful approximation technique over the years. 
See \cite{BW98} for the mathematical validity of the SL method in some cases.
In the context of a basic illustrative example, 
it consists in replacing an equation of the form
$
g(Y_t,\dot Y_t, \ddot Y_t) = X_t
$ 
with an equivalent form
$
m^\star \ddot Y_t + c^\star \dot Y_t + k^\star Y_t = X_t
$
where $Y,\dot Y, \ddot Y$ are, respectively, the displacement response, the velocity and the acceleration, $g$ is a (possibly non-linear) function and $X$ is an input forcing.
The parameters (that can depend on the time) $(m^\star,c^\star,k^\star)$ are determined by minimizing the error 
$\epsilon(c,m,k) \triangleq g(Y_t,\dot Y_t, \ddot Y_t)- m \ddot Y_t - c \dot Y_t - k Y_t$ as follows: $\| \epsilon(c^\star,m^\star,k^\star) \| = \min_{c,m,k \in \mathbb{R}} \| \epsilon(c,m,k) \|$
where $\| .\|$ is a convenient norm.
Further details can be found in \cite{RS03}, with an introduction of SL for simple systems in Chapter 5 and analysis of systems with multiple degrees of freedom in Chapter 6.
The method can be employed to deal with complex systems (including hysteretic elements) having many degrees of freedom and a broad class of excitation (even non-stationary), SL usually gives reasonably good results when non-linear effects are present.
\subsection{Stochastic averaging}
The method of stochastic averaging, originally introduced by \cite{S63}, provides approximate solutions to problems involving the vibration response of lightly damped systems to broad-band random excitation. See \cite{K66,K68,P73,P74} for mathematical rigorous foundations.
For illustration, we consider one-dimensional variables and the presentation we use is inspired by \cite{C09}. Consider a system of the form 
$$
\begin{cases}
& \textup{d} X^\epsilon_t = \epsilon^2 f(X^\epsilon_t,Y^\epsilon_t) \textup{d} t + \epsilon g(X^\epsilon_t,Y^\epsilon_t) \textup{d} W_t, \: X^\epsilon_0 = x,\\ 
& \textup{d} Y^\epsilon_t = h(X^\epsilon_t,Y^\epsilon_t) \textup{d} t + \varphi(X^\epsilon_t,Y^\epsilon_t) \textup{d} W_t, \: Y^\epsilon_0 = y
\end{cases}
$$ 
for $0 < \epsilon \ll 1$ and $f,g,h,\varphi$ functions with appropriate conditions.
Here, $W$ is a Wiener process.
The variables $X^\epsilon_t$ and $Y^\epsilon_t$ are, respectively, called \textit{slow} and \textit{fast} components.
Engineers are interested in the behavior of $X^\epsilon_t$ on intervals of order $\epsilon^{-1}$ since on those time scales the most significant changes occur.
The \textit{averaging principle} states that a good approximation of the slow component can be obtained by averaging over the fast components in the sense that the trajectory $X^\epsilon_t$ is replaced by the solution $\bar{X}$ of 
$$
\textup{d} \bar{X}_t = \epsilon^2 \bar{f}(\bar{X}_t) \textup{d} t + \epsilon \bar{g}(\bar{X}_t) \textup{d} W_t, \: \bar{X}^\epsilon_0 = x,
$$
where 
$\bar{f}(x) = \lim \limits_{T \to \infty} \frac{1}{T} \mathbb{E} \int_0^T f(x,Y_t) \textup{d} t$
and
$\bar{g}^2(x) =\lim \limits_{T \to \infty} \frac{1}{T} \mathbb{E} \int_0^T g^2(x,Y_t) \textup{d} t$
and 
$$
d Y (t) = h(x,Y_t) \textup{d} t + \varphi(x,Y_t) \textup{d} W_t, \: Y^\epsilon_0 = y.
$$
In the dynamics of $Y_t$, $x$ is a frozen parameter.
The stochastic averaging method is particularly useful in dealing with situations with fairly light non-linear damping.
Details on the application of this method to systems of the form $\ddot X + \epsilon^2 \psi(X,\dot X) + k X = \epsilon \dot W$ are presented in \cite{RS86}.
To apply the averaging principle the joint response $(X,\dot X)$ needs to be transformed in 
$X_t = a_t \cos ( \omega_0 t + \varphi_t)$ and $\dot X_t = - a_t \omega \cos (\omega_0 t + \varphi_t)$.
This leads to a stochastic differential equation for $(a_t,\varphi_t)$ which can be seen as the slow component in the equation above.
Details can be found in Section 3 of \cite{RS86}. The stochastic averaging method can cope with systems having hysteretic features such as the elasto-plastic and impact problems.

\subsection{Numerical path integration}
The path integration (PI) method, introduced in \cite{NM96}, is a step-by-step calculation of the joint probability density function (PDF) of a set of state space variables describing a white noise excited nonlinear dynamical system. Precisely, the PDF is computed at a given time by applying the Chapman-Kolmogorov equation when the response PDF and the transition PDF are known at an earlier time.
The PI method is an efficient approximation for solving the Fokker-Planck equation and for providing the stationary response of the underlying dynamical system.
The method exploits short time transition probability density functions which can be sharply approximated by explicit formulae (based on a local Gaussian behavior) even for hysteretic systems such as elasto-plastic or impact problems.

\subsection{Comments}
The statistical linearization and stochastic averaging methods rely on an approximation of the underlying stochastic process, where in our method no such approximation is needed.
However numerical path integration methods deal with the original process. Our approach relies heavily on the PDEs related to the Feynman-Kac (FK) formula and the backward Kolmogorov equations whereas the PI method uses the Chapman-Kolmogorov equation which ultimately lead to the Forward Kolmogorov equation. Hence PI and FK methods can be seen as being the dual of each other. 
In contrast with the PI method, our method extends naturally (by adapting the infinitesimal generator) to solving a broad range of problems. These include optimal stopping and stochastic optimal control problems via free boundary value problem and HJB equations with non-standard boundary condition, see \cite{LLMWZ19}.
Our approach can be combined with existing methodologies to investigate a number of practically important quantities.
For example, following the steps of \cite{SKS11} and \cite{BB19} it can be used to study the power spectrum density for stationary processes of the form
$$
S(\omega) \triangleq \lim_{T \to \infty} \frac{1}{T} \mathbb{E} \left ( \left | \int_0^T  g(X_t,Y_t) e^{-\mathbf{i} \omega t} \textup{d} t  \right |^2 \right  ),
$$
where $X,Y$ is the response variable of one of the systems studied in this paper (elasto-plastic or impact system) and the notation $|.|$ stands for the modulus of a complex number. 
Some details are given in the conclusions and perspectives section.
}

\section{The Partial Differential Equations for $A, B ,C$ and $A',B', C'$}
\noindent It is possible to relate the functions $A$ and $A'$ to the solution of \textit{parabolic problems},
the functions $B$ and $B'$ to the solution of \textit{elliptic problems} and the functions $C$ and $C'$ to the solution of \textit{Poisson problems}. In the first part of this section, we present these problems in the case of $\mathbf{H} = 0$. Then, in the second part, we provide a formal presentation of the corresponding problems for the two non-smooth problems targeted in this paper.\\

\subsection*{Notation}
\noindent For $T>0$ and a domain $\Omega$ of $\mathbb{R}^2$, we use the notation $C^\star(\Omega \times [0,T])$ for the set of continuous functions on $\Omega \times [0,T]$ that are $C^1$-regular with respect to $x$, $C^2$-regular with respect to $y$ and $C^1$-regular with respect to $t$. We use the notation $C^\star (\Omega)$ for the set of continuous functions on $\Omega$ that are $C^1$-regular with respect to $x$ and $C^2$-regular with respect to $y$. 
We use the generic notation $\Gamma$ for 
\[
\Gamma_{T}^\lambda(f,g) 
\triangleq e^{-\lambda T} f(X_T,Y_T)  
+ \int_0^T e^{-\lambda \tau} g(X_\tau,Y_\tau) \d \tau
\]
because it helps to write the quantities $A,B,C$ and $A',B',C'$ in a compact form.\\

\noindent \lmxxx{We use the following notations below  
\begin{equation}
\label{A0}
\tag{$\mathcal{A}_0(\lambda, \phi)$}
\phi \in C^\star (\mathbb{R}^2 \times [0,\infty])
\quad
\mbox{and}
\quad
\mathbb{E} \left ( \int_0^t e^{-2 \lambda \tau} \left | \frac{\partial \phi }{\partial y} \right |^2(X_\tau,Y_\tau, \tau) \d \tau \right ) < \infty, \forall t. 
\end{equation}
\begin{equation}
\label{A1}
\phi \in C^\star (\mathbb{R}^2)
\quad
\mbox{and}
\quad
\tag{$\mathcal{A}_1(\lambda,\phi)$}
\lim \limits_{T \to \infty } \exp(-\lambda T) \mathbb{E} \phi(X_T,Y_T) = 0. 
\end{equation}
\begin{equation}
\label{A2}
\phi, \psi \in C^\star (\mathbb{R}^2)
\quad
\mbox{and}
\quad
\tag{$\mathcal{A}_2(\lambda,\phi)$}
\int_0^\infty \exp(-\lambda \tau ) \mathbb{E} | \phi | (X_\tau,Y_\tau) \d \tau < \infty.
\end{equation}
\begin{equation}
\label{A3}
\phi, \psi \in C^\star (\mathbb{R}^2)
\quad
\mbox{and}
\quad
\tag{$\mathcal{A}_3(\lambda, \mu, \phi, \psi)$}
\lim \limits_{T \to \infty } \exp(- \lambda T) \mathbb{E} \left ( \phi (X_T,Y_T) \int_0^T \exp(-\mu \tau ) \psi(X_\tau,Y_\tau) \d \tau \right ) = 0.
\end{equation}
}
\subsection{The case $\mathbf{H} = 0$}
\noindent  \lmxxx{The proofs of connecting the PDEs below and stochastic processes can be found in the Appendix for the reader's convenience.}
\textcolor{black}{We use the notation 
\[
L  \triangleq \frac{1}{2} \frac{\partial^2 }{\partial y^2} 
+ F \frac{\partial }{\partial x} + G \frac{\partial }{\partial y}.
\]}
\subsubsection*{Backward-in-time parabolic problems}
\label{THMA}
Consider functions $u \in C^\star (\mathbb{R}^2 \times [0,T]), v \in C^\star (\mathbb{R}^2 \times [0,T+h]), w \in C^\star (\mathbb{R}^2 \times [0,T])$. 
Let $\lambda, \mu \geq 0.$ 
Assume that $u,v,w$ satisfy $\mathcal{A}_0(\lambda, u)$, $\mathcal{A}_0(\mu, v)$, $\mathcal{A}_0(\lambda + \mu, w)$, respectively, and that 
\begin{equation}
\label{cu}
\frac{\partial u}{\partial t} +  Lu - \lambda u = -g, \quad u(T) = f \quad \mbox{in} \quad \mathbb{R}^2,
\quad 
\frac{\partial v}{\partial t} +  Lv - \mu v = -\psi, \quad v(T+h) = \varphi \quad \mbox{in} \quad \mathbb{R}^2,
\end{equation}
and
\begin{equation}
\label{cw}
\frac{\partial w}{\partial t} +  Lw - (\lambda + \mu) w = -\frac{\partial u}{\partial y} \frac{\partial v}{\partial y},
\quad w(T) = 0 \quad \mbox{in} \quad \mathbb{R}^2.
\end{equation}
Then we have 
\[
A = u(x,y,0) \quad \mbox{and} \quad A' = (uv + w)(x,y,0).
\] 
\noindent As a corollary,
\[
\textup{Var}(\Gamma_{T}^0(f,g)) = w(x,y,0),
\]
where $w$ solves the problem \eqref{cw} with $v = u$ and $h=0$. 
\subsubsection*{Degenerate elliptic problems}
\label{THMB}
\lmxxx{Let $\lambda, \mu \geq 0$ be two numbers and $g, \psi$ be two functions such that 
$$
\int_0^\infty \exp(-\lambda \tau ) \mathbb{E} |g|(X_\tau,Y_\tau) \d \tau < \infty, 
\quad  \int_0^\infty \exp(-\mu \tau ) \mathbb{E} |\psi|(X_\tau,Y_\tau) \d \tau < \infty.
$$ 
Consider functions $u_\lambda, v_\mu , w_{\lambda+\mu} \in C^\star (\mathbb{R}^2)$. 
Assume that $u_\lambda,v_\mu,w_{\lambda+\mu}$ satisfy technical conditions 
$\mathcal{A}_1(\lambda,u_\lambda), \mathcal{A}_1(\mu,v_\mu)$, 
$\mathcal{A}_1(\lambda+\mu,w_{\lambda+\mu})$, 
$\mathcal{A}_1(\lambda, \mu, \frac{\partial u_\lambda}{\partial y} \frac{\partial v_\mu}{\partial y})$, 
$\mathcal{A}_1(\lambda+\mu, u_\lambda v_\mu)$
$\mathcal{A}_3(\lambda, \mu, u_\lambda, \psi)$,
$\mathcal{A}_3(\mu, \lambda, v_\mu, g)$}
and
\begin{equation}
\label{cuvw}
\begin{dcases}
& Lu_\lambda - \lambda u_\lambda = -g \quad \mbox{in} \quad \mathbb{R}^2,\\
& Lv_\mu - \mu v_\mu = -\psi \quad \mbox{in} \quad \mathbb{R}^2,\\
& Lw_{\mu + \lambda} - (\mu + \lambda) w_{\lambda + \mu} = -\frac{\partial u_\lambda}{\partial y} \frac{\partial v_\mu}{\partial y}
\quad \mbox{in} \quad \mathbb{R}^2.
\end{dcases}
\end{equation}
Then we have 
\[
B = u_\lambda(x,y) \quad \mbox{and} \quad B' = (u_\lambda v_\mu + w_{\lambda + \mu})(x,y).
\] 

\subsubsection*{Long time behavior and degenerate Poisson problems}
\label{THMC}
Assume that $(X_t,Y_t)$ has a unique invariant probability measure $\nu$ on $\mathbb{R}^2$ in the sense that  
for any continuous function $f$ satisfying $\nu(|f|) < \infty$, we have $\mathbb{E} f(X_t,Y_t) = \nu(f)$ provided that $(X_0,Y_0)$ is distributed according to $\nu$.\\

\noindent Consider functions $U, V \in C^\star (\mathbb{R}^2)$. 
Assume that $g$ and $\psi$ satisfy $\nu(|g|) < \infty, \nu(|\psi|) < \infty$ and $U, V$ satisfy 
\[
\nu \left ( \left  | \frac{\partial U}{\partial y} \frac{\partial V}{\partial y} \right |  \right ) < \infty
\]
together with
\begin{equation}
\label{cuvw}
L U = \nu(g) - g \quad \mbox{in} \quad \mathbb{R}^2, \quad  LV  = \nu(\psi) - \psi \quad \mbox{in} \quad \mathbb{R}^2.
\end{equation}
Then we have 
\[
C = \nu (g) \quad \mbox{and} \quad C' = \nu \left ( \frac{\partial U}{\partial y} \frac{\partial V}{\partial y} \right ).
\] 

\subsection{Non-standard PDEs of SVIs: elasto-plastic and obstacle problems}
\noindent In a similar manner to the results above for the case $\mathbf{H}=0$, here we present non-standard PDEs related to SVIs modeling the elasto-plastic \eqref{pathdependent_ep} and obstacle \eqref{pathdependent_op} problems. 
First, for each problem, we give a description of the infinitesimal generator of the corresponding process.
Then, we only present the non-standard backward-in-time parabolic problems in analogy to Section \ref{THMA}. 
For the non-standard elliptic degenerate and Poisson problems, the idea remains basically the same as what is presented above.\\ 

\noindent Without loss of generality, we can assume here that $P_Y = P_O = 1$, define
\[
D \triangleq (-1,1) \times (-\infty, \infty), \quad D_{T} \triangleq (-1,1) \times (-\infty, \infty) \times (0,T)
\]
and
\[
D_{\pm} \triangleq \{ \pm 1 \} \times (-\infty, \infty), \quad D_{T}^\pm \triangleq \{ \pm 1 \} \times (-\infty, \infty) \times (0,T).
\] 
\noindent 
In the elasto-plastic problem, using Ito's lemma, the generator of $(Z_t,Y_t)$ is defined on any function $\phi \in C^\star(\bar{D})$ and satisfies  
\[
\lim_{t \to 0} \frac{\mathbb{E} \phi(Z_t,Y_t) - \phi(z,y)}{t} =
\left \{
\begin{array}{rcl}
L \phi \triangleq \dfrac{1}{2} \dfrac{\partial^2 \phi}{\partial y^2} - (c_0 y + kz) \dfrac{\partial \phi}{\partial y} + y \dfrac{\partial \phi}{\partial z}, & \mbox{ if } & |z| <1,\\[2mm] 
L_+ \phi \triangleq \dfrac{1}{2} \dfrac{\partial^2 \phi}{\partial y^2} - (c_0 y + k) \dfrac{\partial \phi}{\partial y} + \min(0,y) \dfrac{\partial \phi}{\partial z}, & \mbox{ if } & z=1,\\[2mm] 
L_- \phi \triangleq \dfrac{1}{2} \dfrac{\partial^2 \phi}{\partial y^2} - (c_0 y - k) \dfrac{\partial \phi}{\partial y} + \max(0,y) \dfrac{\partial \phi}{\partial z}, & \mbox{ if } & z=-1. 
\end{array}
\right.
\]
In the obstacle problem, (this is formal) the generator of $(X_t,Y_t)$ is defined on any function $\phi \in C^\star(\bar{D})$ such that 
\[
\phi(\pm 1, y) = \phi(\pm 1, -e y) \: \mbox{in} \: \pm y >0
\]
satisfies  
\[
\lim_{t \to 0} \frac{\mathbb{E} \phi(X_t,Y_t) - \phi(x,y)}{t} =
\dfrac{1}{2} \dfrac{\partial^2 \phi}{\partial y^2} - (c_0 y + kx) \dfrac{\partial \phi}{\partial y} + y \dfrac{\partial \phi}{\partial x}, \mbox{ if }  |x| <1.
\]

\subsubsection*{Non-standard problems for the elasto-plastic problem}
\label{THM_ep}
Consider functions $u \in C^\star ( D \times [0,T]), v \in C^\star (D \times [0,T+h]), w \in C^\star (D \times [0,T])$. 
Assume that $u,v,w$ satisfy $\mathcal{A}_0(\lambda, u)$, $\mathcal{A}_0(\mu, v)$, $\mathcal{A}_0(\lambda + \mu, w)$, respectively, (w.r.t to the solution of \eqref{svi_ep})
\begin{equation}
\label{uep}
\frac{\partial u}{\partial t} +  Lu - \lambda u = -g \: \mbox{in} \: D_T, 
\quad \frac{\partial u}{\partial t} +  L_\pm u - \lambda u = -g \: \mbox{in} \: D_T^\pm,  
\quad u(T) = f  \: \mbox{in} \: D.
\end{equation}
\begin{equation}
\label{vep}
\frac{\partial v}{\partial t} +  Lv - \mu v = -\psi \: \mbox{in} \: D_{T+h}, 
\quad \frac{\partial v}{\partial t} +  L_\pm v - \mu v = -\psi \: \mbox{in} \: D_{T+h}^\pm,  
\quad v(T+h) = \phi  \: \mbox{in} \: D.
\end{equation}
\begin{equation}
\label{wep}
\frac{\partial w}{\partial t} +  Lw - (\lambda + \mu) w = -\frac{\partial u}{\partial y} \frac{\partial v}{\partial y}  \: \mbox{in} \: D_T, 
\quad \frac{\partial w}{\partial t} +  L_\pm w - (\lambda + \mu) w = - \frac{\partial u}{\partial y} \frac{\partial v}{\partial y} \: \mbox{in} \: D_T^\pm,  
\quad w(T) = 0  \: \mbox{in} \: D
\end{equation}
then the corresponding $A$ and $A'$ satisfy
\[
A = u(x,y,0) \quad \mbox{and} \quad A' = (uv + w)(x,y,0).
\] 
\subsubsection*{Non-standard problems for the obstacle problem}
\label{THM_op}
Consider functions $u \in C^\star (D \times [0,T]), v \in C^\star (D \times [0,T+h]), w \in C^\star (D \times [0,T])$. 
Assume that $u,v,w$ satisfy $\mathcal{A}_0(\lambda, u)$, $\mathcal{A}_0(\mu, v)$, $\mathcal{A}_0(\lambda + \mu, w)$, respectively, (w.r.t to the solution of \eqref{svi_obstacle})
\begin{equation}
\label{uop}
\frac{\partial u}{\partial t} +  Lu - \lambda u = -g \: \mbox{in} \: D_T, 
\quad u(\pm 1, y) = u(\pm 1, -e y) \: \mbox{in} \: D_T^\pm,  
\quad u(T) = f  \: \mbox{in} \: D.
\end{equation}
\begin{equation}
\label{vop}
\frac{\partial v}{\partial t} +  Lv - \mu v = -\psi \: \mbox{in} \: D_{T+h}, 
\quad v(\pm 1, y) = v(\pm 1, -e y) \: \mbox{in} \: D_{T+h}^\pm,
\quad v(T+h) = \varphi \: \mbox{in} \: D.
\end{equation}
\begin{equation}
\label{wop}
\frac{\partial w}{\partial t} +  Lw - (\lambda + \mu) w = - \frac{\partial u}{\partial y} \frac{\partial v}{\partial y}  \: \mbox{in} \: D_T, 
\quad w(\pm 1, y) = w(\pm 1, -e y) \: \mbox{in} \: D_T^\pm,  
\quad w(T) = 0  \: \mbox{in} \: D
\end{equation}
then the corresponding $A$ and $A'$ satisfy
\[
A = u(x,y,0) \quad \mbox{and} \quad A' = (uv + w)(x,y,0).
\] 
\lmxxx{
\section{Asymptotic formula for D}
\noindent 
For every function $f$ on $\mathbb{R}^2$, define
$
\Delta_t^{EP}(f) \triangleq \int_0^t f(Z_s,Y_s) \d s
\: \mbox{ and } \:  
\Delta_t^{OP}(f) \triangleq \int_0^t f(X_s,Y_s) \d s 
$
where $(Z_s,Y_s)$ solves the elasto-plastic problem and $(X_s,Y_s)$ solves the obstacle problem. 
The case \eqref{epp_crossings} corresponds to the case $f(z,y) = y$. 
The case \eqref{op_crossings} corresponds to the case $f(x,y) = x$.
We drop the superscript $(.)^{EP/OP}$ for convenience. In both cases, the process $t \mapsto \Delta_t(f)$ has continuous trajectories and admits an asymptotic (in time) variance of the form $\gamma_f^2 t, \: \gamma_f^2 >0$.
Thus, with respect to an appropriate scaling in time and space, it is reasonable to expect that it behaves like a Brownian motion in the sense that, taking $(b,T) = (\sqrt{p}\, \beta, p \Theta), \: p$ large, the following approximation can be made
\begin{equation}
\label{equivalent}
\lim_{p \to \infty} W_f (\sqrt{p}\, \beta, p \Theta)
= 1- \frac 4 \pi \sum_{k=0}^{+\infty} \frac {(-1)^k} {2k+1} \exp\left( -\frac{(2k+1)^2\pi^2 \Theta}{8 \beta^2} \gamma_f^2 \right),
\end{equation}
where
\[
\gamma_f^2 \triangleq \lim_{t \to \infty} \frac{\sigma^2(\Delta_t(f))}{t}
\quad
\mbox{and}
\quad
W_f(b,T) \triangleq \mathbb{P} \left ( \max_{0 \leq t \leq T} \left | \Delta_t(f) \right | \geq b \right ).
\]
For an elasto-perfectly-plastic oscillator excited by white noise, this Brownian behavior was already proposed in \cite{BL89} as a heuristic argument 
and, recently, a functional central limit theorem (FCLT) has been given in \cite{FLM18}. 
In the case of a penalization of variational inequalities, a FCLT has been obtained in a general framework including penalization of \eqref{svi_ep} and \eqref{svi_obstacle} (with purely elastic impacts, $e = 1$) \cite{LM18}. The computation of a quantity of type $D$ relies on the computation of quantity of type $C'$ (shown in Section~4). 
See Figure~\ref{fig:D} for illustration.
\begin{figure}[h!]
\centering
   \begin{subfigure}[t]{0.4\linewidth}
          \begin{tikzpicture}[scale=0.66]
        \begin{axis}[legend style={at={(1,0)},anchor=north west}, compat=1.3,
          xmin=0, xmax=3,ymin=0,ymax=1,
          xlabel= {$T$},
          ylabel= {$$},
          tick label style={/pgf/number format/fixed},
          legend style={at={(1,0)},anchor=south east},
          legend cell align=left]
          \addplot[thick,solid,color=red,mark=none] table [x index=0, y index=1]{./explicit_tmax5_data.txt};
          \addplot[thick, dotted,color=black,mark=none] table [x index=0, y index=1]{./WMC_p1.0E+01_dt1.0E-03_MC_1.0E+03_data.txt};
          \addplot[thick, dashdotted,color=black,mark=none] table [x index=0, y index=1]{./WMC_p1.0E+02_dt1.0E-03_MC_1.0E+03_data.txt};
          \addplot[thick, dashed,color=black,mark=none] table [x index=0, y index=1]{./WMC_p1.0E+03_dt1.0E-03_MC_1.0E+03_data.txt};
         \legend{RHS of \eqref{equivalent}, $p=10$, $p=100$, $p=1000$}
        \end{axis}
        \end{tikzpicture}
        \caption{\textbf{Elasto-plastic problem.} 
        We consider $P_Y = 0.25$.
        The dotted and dashed lines represent the Monte Carlo approximation of $W_y(\sqrt{p}\beta, p\Theta)$, for $3$ different values of $p$, 
        while the solid line represent the explicit formula shown in the right hand side of \eqref{equivalent} with $f(z,y) = y$. 
        	\label{fig:D-C}}
\end{subfigure}
 \quad
   \begin{subfigure}[t]{0.4\linewidth}
          \begin{tikzpicture}[scale=0.66]
        \begin{axis}[legend style={at={(1,0)},anchor=north west}, compat=1.3,
          xmin=0, xmax=3,ymin=0,ymax=1,
          xlabel= {$T$},
          ylabel= {$$},
          tick label style={/pgf/number format/fixed},
          legend style={at={(1,0)},anchor=south east},
          legend cell align=left]
  \addplot[thick, solid,color=red,mark=none] table [x index=0, y index=1]{./obstacle_explicit_tmax5_data.txt};
  \addplot[thick, dotted,color=black,mark=none] table [x index=0, y index=1]{./obstacle_WMC_p1.0E+01_dt1.0E-04_MC_1.0E+03_data.txt};     
  \addplot[thick, dashdotted,color=black,mark=none] table [x index=0, y index=1]{./obstacle_WMC_p1.0E+02_dt1.0E-04_MC_1.0E+03_data.txt};   
  \addplot[thick, dashed,color=black,mark=none] table [x index=0, y index=1]{./obstacle_WMC_p1.0E+03_dt1.0E-04_MC_1.0E+03_data.txt};       
       \legend{RHS of \eqref{equivalent}, $p=10$, $p=100$, $p=1000$}
        \end{axis}
        \end{tikzpicture}
        \caption{\textbf{Obstacle problem.} 
        We consider $P_0 = 1$. The dotted and dashed lines represent the Monte Carlo approximation of $W_x(\sqrt{p}\beta, p\Theta)$, for $3$ different values of $p$, 
        while the solid line represent the explicit formula shown in the right hand side of \eqref{equivalent} with $f(x,y) = x$.   
        	\label{fig:D-D}}
\end{subfigure}
\caption{\label{fig:D} Numerical results for $\eqref{epp_crossings}$ and $\eqref{op_crossings}$. In both cases, $\gamma_f^2$ is computed using the PDE method.
For the Monte Carlo simulation, we used $10^5$ samples with a time step of $\delta t = 10^{-4}$. Here $\beta=0.5, \Theta=3$, thus $b = \sqrt{p} \beta$ and $T = p \Theta$.}
\end{figure}
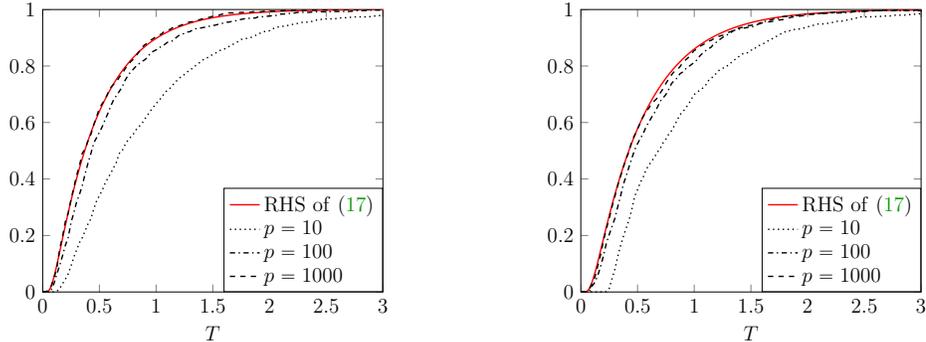
}
\subsection*{Computational savings}
\lmxxx{In order to compute the results shown in Figures 2a and b, we used MATLAB codes (running on Intel Xeon E5 2.3Ghz processes) for both our proposed method and the Monte-Carlo method. Our method typically took around one minute, whereas the Monte-Carlo method took around 5 hours to achieve similar accuracy.}

\section{Numerical computation of $A, B, C$ and $A',B', C'$ }
\noindent In this section, we first explain the probabilistic numerical approach for SVIs of the elasto-plastic \eqref{svi_ep} and obstacle \eqref{svi_obstacle} problems. Then, we present a numerical approach for solving \lmxxx{the PDEs with non-standard boundary condition related to the FKf}. The method is first presented and applied to the elasto-plastic problem \eqref{uep}, \eqref{vep} and \eqref{wep}.
Then, it is applied to the obstacle problem, \eqref{uop}, \eqref{vop}
and \eqref{wop}.

\subsection{Discretization of SVIs and Monte Carlo approach}
\noindent For the probabilistic numerical scheme of $\{ (Z_t,Y_t), \:
t \geq 0 \}$ of \eqref{svi_ep}, we use a time step $\delta t >0$.
Then we construct random variables $\{ (Z_n, Y_n), 1 \leq n \leq
N_{\delta t}\}$ and a partition of $[0,T]$, $0 = t_0 < t_1 < t_2 < \ldots < t_{N_{\delta t}} = T$.
Here, for each $1 \leq n \leq N_{\delta t}, (Z_n, Y_n)$ is an
approximation of $(Z_{t_n}, Y_{t_n})$. We refer to \cite{BMPT09} for a detailed
presentation of the algorithm. To be concise, we only explain the
computation of quantities of type $A$. For any choice of
well-behaved functions $f$ and $g$, we proceed with the following
approximation
\begin{equation}
\label{approxproba}
\mathbb{E} \left ( f(Z_T,Y_T) + \int_0^T g(Z_\tau,Y_\tau) \d \tau \right ) 
\approx  \frac{1}{M} \sum_{m=1}^M 
\left (
f (Z_{N_{\delta t}}^m,Y_{N_{\delta t}}^m)
+
\sum_{n=0}^{N_{\delta t} - 1} g(Z_n^m,Y_n^m) (t_{n+1}-t_n)
\right )
\end{equation}
where $\{ (Z^m, Y^m), m = 1, \ldots, M \}$ is an i.i.d.\ sequence of
trajectories produced by the algorithm.\\ 

\noindent Details of the probabilistic numerical scheme to compute $\{
(X_t,Y_t), \: t \geq 0 \}$ in \eqref{svi_obstacle} are given in
Appendix~\ref{algo_svi_obstacle}.  Then we use a similar approximation to that given by \eqref{approxproba} in terms of $(X,Y)$ of \eqref{svi_obstacle}. 
\subsection{Discretization of PDE problems}
\label{ssec:num1}
\noindent 
To numerically approximate the solutions of \eqref{uep}, \eqref{vep} and \eqref{wep}, we use a finite difference scheme.
We truncate the unbounded domain $D$ to obtain $D_Y \triangleq (-1,1) \times (-Y,Y)$, where $Y$ is chosen sufficiently large that the probability of finding the underlying process outside $D_Y$ is negligible. We apply a homogeneous Neumann boundary condition at $y = \pm Y$ . We consider a two-dimensional rectangular finite difference grid, 
\[
\mathcal{G} \triangleq \left \{ (z_i,y_j) \triangleq (-1 + (i-1)  \delta z,-Y + (j-1) \delta y) \right \}_{1 \leq i \leq I, 1 \leq j \leq J},
\]
where $\delta z \triangleq \frac{2}{I-1}, \delta y \triangleq \frac{2Y}{J-1}$. Here, $I,J$ are odd integers of the form $2\tilde{I}+1,2\tilde{J}+1$. 
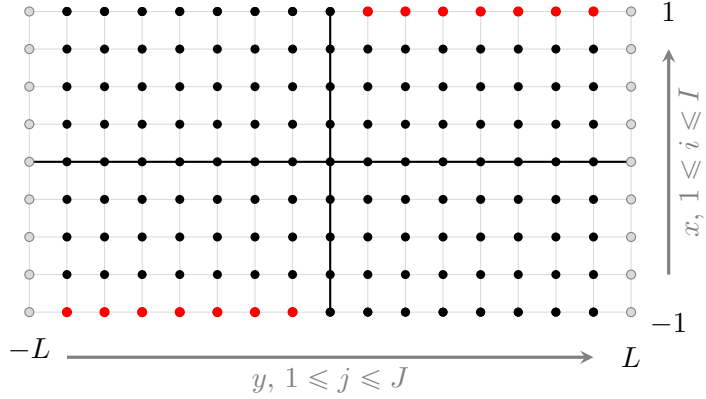
\begin{figure}
\begin{tikzpicture}[scale=0.5,rotate=90]
\draw[step=1cm, gray!30!white, very thin] (-4,-8) grid (4,8);
\draw[thick] (-4,0) -- (4,0);
\draw[thick] (0,-8) -- (0,8);
\foreach \x in {-4,...,4}
   \foreach \y in {-7,...,7}
   \fill[black] (\x,\y) circle (0.12cm);
\foreach \x in {-4}
   \foreach \y in {-7,...,0}
   \fill[black] (\x,\y) circle (0.12cm);
\foreach \x in {4}
   \foreach \y in {0,...,7}
   \fill[black] (\x,\y) circle (0.12cm);
\foreach \x in {-4,-3,...,4}
   \foreach \y in {-8,8}{
   \fill[gray!30!white] (\x,\y) circle (0.12cm);
   \draw[gray] (\x,\y) circle (0.12cm);
}
\foreach \x in {4}
   \foreach \y in {-7,...,-1}{
   \fill[red] (\x,\y) circle (0.12cm);
   \draw[red] (\x,\y) circle (0.12cm);
}
\foreach \x in {-4}
   \foreach \y in {1,...,7}{
   \fill[red] (\x,\y) circle (0.12cm);
   \draw[red] (\x,\y) circle (0.12cm);
}
\node at (-4.3,-9) {\small $-1$};
\node at (4,-9) {\small $1$};
\node at (-5.2,-8) {\small $L$};
\node at (-5,8) {\small $-L$};
\draw[->,>=stealth,very thick,gray] (-5.2,7) -> (-5.2,-7);
\node at (-5.8,0) {\small\textcolor{gray}{ $y$, $1\le j\le J$}};
\draw[->,>=stealth,very thick,gray] (-3,-9.) -> (3,-9.);
\node[rotate=90] at (0,-9.6) {\small\textcolor{gray}{ $x$, $1\le i\le I$}};
\end{tikzpicture}
\caption{Discretization of $D$. At black points, the discretized
  equation is satisfied. At grey points, homogeneous Neumann boundary
 conditions are used, and at red points non-standard boundary
 conditions are employed.\label{fig:grid}}
\end{figure}
The total number of nodes in $\mathcal{G}$ is $N = IJ$. The numerical
approximations of $u(z_i,y_j,t_n)$, $v(z_i,y_j,t_n)$ and
$w(z_i,y_j,t_n)$ are denoted by ${u}_{i,j}^n$,
${v}_{i,j}^n$ and ${w}_{i,j}^n$ and the corresponding vectors
collecting all the unknowns are $\boldsymbol{u}^n$,
$\boldsymbol{v}^n$ and $\boldsymbol{w}^n$. We use the notation $f_{i,j}, g_{i,j}, \psi_{i,j}, \varphi_{i,j}$ for $f(x_i,y_j), g(x_i,y_j), \psi(x_i,y_j)$ 
and $\varphi(x_i,y_j)$ and the corresponding vectors are 
$\boldsymbol{f}$, $\boldsymbol{g}$, $\boldsymbol{\varphi}$, $\boldsymbol{\psi}$. Here, $t_n \triangleq n \delta t$ discretizes the time and $N_T \delta t = T$, $N_{T+h} \delta t = T+h$.
Using an implicit Euler method to discretize in time together with finite differences in space, the first conditions in \eqref{uep}, \eqref{vep} and \eqref{wep} at the black points in Figure \ref{fig:grid} result in
\[
\begin{dcases}
& {u}_{i,j}^0 = f_{i,j} \quad \mbox{and} \quad \frac{{u}_{i,j}^{n+1}-{u}_{i,j}^{n}}{\delta t} - \left ( L \boldsymbol{u}^{n+1} \right )_{i,j} + \lambda {u}_{i,j}^{n+1} = g_{i,j}, \quad 1 \leq n \leq N_T-1,\\ 
& {v}_{i,j}^0 = \varphi_{i,j} \quad \mbox{and} \quad \frac{{v}_{i,j}^{n+1}-{v}_{i,j}^{n}}{\delta t} - \left ( L \boldsymbol{v}^{n+1} \right )_{i,j} + \mu {v}_{i,j}^{n+1} = \psi_{i,j}, \quad 1 \leq n \leq N_{T+h}-1,\\ 
& {w}_{i,j}^0 = 0 \quad \mbox{and} \quad
\frac{{w}_{i,j}^{n+1}-{w}_{i,j}^{n}}{\delta t} - \left ( L \boldsymbol{w}^{n+1} \right )_{i,j} = 
\left ( D_y \boldsymbol{u}^n \right )_{i,j} \left ( D_y \boldsymbol{v}^n \right )_{i,j}, \quad 1 \leq n \leq N_T-1,\\   
\end{dcases}
\]
where $\left ( D_y \boldsymbol{u} \right )_{i,j}$ and $(L \boldsymbol{u})_{i,j}$ are centered finite differences,
\[
\left ( D_y \boldsymbol{u} \right )_{i,j} \triangleq \frac{{u}_{i,j+1}-{u}_{i,j-1}}{2 \delta y}
\]
and first order upwind differences
\begin{equation}
\label{discret_mL}
-(L \boldsymbol{u})_{i,j} \triangleq  -(L_y \boldsymbol{u})_{i,j} - \max(0, y_j) \left ( \frac{u_{i+1,j} - u_{i,j}}{\delta z} \right ) 
- \min(0, y_j) \left ( \frac{u_{i,j} - u_{i-1,j}}{\delta z} \right )
\end{equation}
with
\[
-(L_y \boldsymbol{u})_{i,j} \triangleq - \dfrac{1}{2} \left (\dfrac{u_{i,j+1} - 2 u_{i,j} + u_{i,j-1}}{\delta y^2} \right )
 -\max(0,b_{i,j}) \left ( \frac{u_{i,j+1} - u_{i,j} }{\delta y} \right ) - \min(0,b_{i,j}) \left ( \frac{u_{i,j} - u_{i,j-1} }{\delta y} \right ).
 \]
The nonstandard  boundary conditions (second condition) in \eqref{uep}, \eqref{vep} and \eqref{wep} at the red points in Figure \ref{fig:grid} are discretized by the same formulae with $L$ replaced by $L_\pm$. Here,  $(L_{\pm} \boldsymbol{u})_{i,j}$ are defined by
\begin{equation}
\label{discret_mLp}
-(L_+ \boldsymbol{u})_{i,j} 
\triangleq
-(L_y \boldsymbol{u})_{i,j}
- \min(0, y_j) \left ( \frac{u_{i,j} - u_{i-1,j}}{\delta z} \right )
\end{equation}
and
\begin{equation}
\label{discret_mLm}
-(L_- \boldsymbol{u})_{i,j} \triangleq  -(L_y \boldsymbol{u})_{i,j} - \max(0, y_j) \left ( \frac{u_{i+1,j} - u_{i,j}}{\delta z} \right ). 
\end{equation}
The Neumann boundary conditions at the points shown in grey results in
\begin{equation}
\label{n3}
(N_2 \boldsymbol{u})_{i,j} = 0,
\quad (N_2 \boldsymbol{v})_{i,j} = 0, 
\quad (N_2 \boldsymbol{w})_{i,j} = 0,
\quad \mbox{where} \quad (N_2 \boldsymbol{u})_{i,j} \triangleq \left \{ \begin{matrix} \dfrac{u_{i,j+1}-u_{i,j}}{\delta y} \quad \text{if } j=1,\\ \\ \dfrac{u_{i,j}-u_{i,j-1}}{\delta y} \quad\text{if } j=J. \end{matrix} \right.  
\end{equation}
This results in the following linear system to be solved in each time step:
\begin{equation}
\label{algorithm}
(I+ \delta t  M) \boldsymbol u^{n+1} = \boldsymbol u^n + \delta t
\boldsymbol g, \quad \boldsymbol u^0 = \boldsymbol f,
\end{equation}
where $M$ is a sparse $N \times N$-matrix that does not depend on $n$.\\

\noindent For the computational results presented in the remainder of this
paper, we use a \texttt{C} code that implements a Monte Carlo (MC)
probabilistic simulation to approximate the solution of
\eqref{svi_ep}.  We use a \texttt{MATLAB} implementation for
the PDE approach \eqref{algorithm}. Here, an LU factorization is used
for solving the linear systems arising from the PDE-approach. For
time-dependent PDE simulations, the LU factors are computed once
upfront and reused throughout the simulation.  Implementations are
available upon request.\\

\subsection{Numerical results for the elasto-plastic problem}
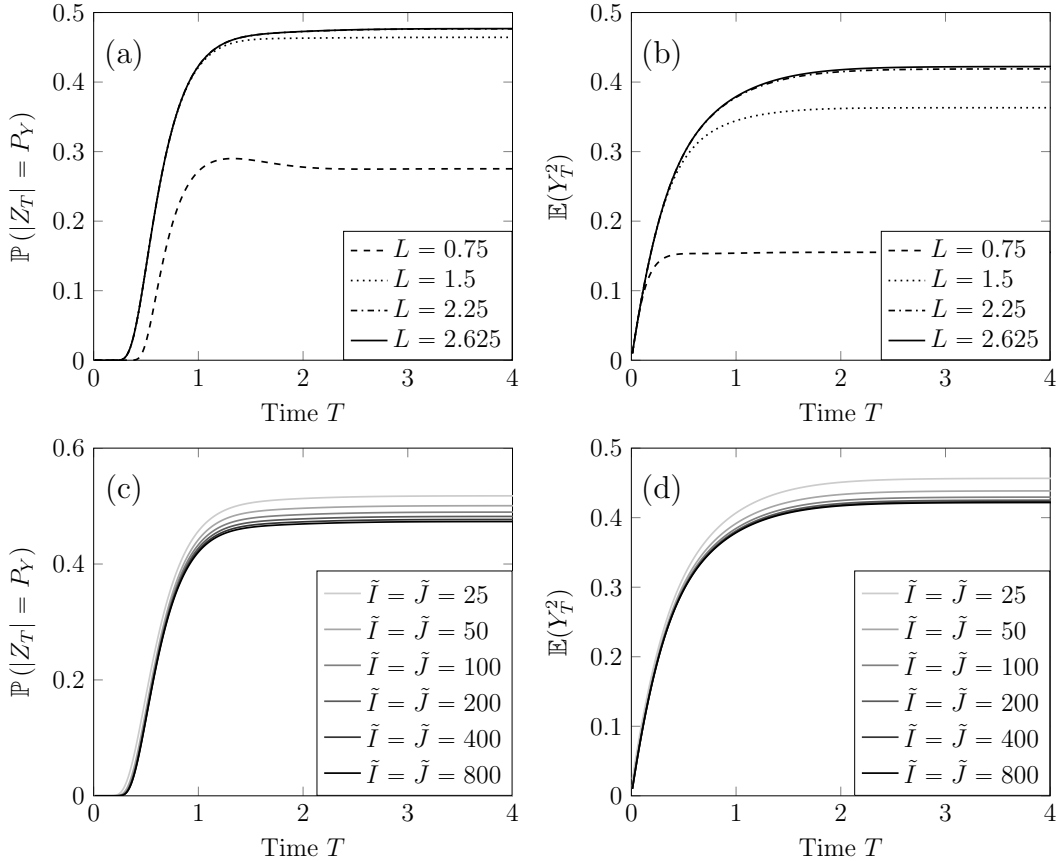
\begin{figure}[th]
\centering
\begin{tikzpicture}[scale=0.8125]
\begin{axis}[legend style={at={(1,0)},anchor=south east}, compat=1.3,
  xmin=0, xmax=4,ymin=0,ymax=0.5,
  xlabel= {Time $T$},
  ylabel= {$\mathbb{P} \left ( |Z_T| = P_Y \right )$},
  tick label style={/pgf/number format/fixed},
  legend cell align=left]
 \addplot[thick, dashed,color=black,mark=none] table [x index=0, y index=1]{./sL_plasticstate_pde100_dt1.0e-03_data.txt};
 \addplot[thick, dotted,color=black,mark=none] table [x index=0, y index=1]{./sL_plasticstate_pde200_dt1.0e-03_data.txt};
 \addplot[thick, dashdotted,color=black,mark=none] table [x index=0, y index=1]{./sL_plasticstate_pde300_dt1.0e-03_data.txt};
 \addplot[thick,color=black,mark=none] table [x index=0, y index=1]{./sL_plasticstate_pde350_dt1.0e-03_data.txt};
 \legend{$L=0.75$,$L=1.5$,$L=2.25$,$L=2.625$}
\end{axis}
\node at (0.5cm,5cm) {\textcolor{black}{(a)}};
\end{tikzpicture}
\begin{tikzpicture}[scale=0.8125]
\begin{axis}[legend style={at={(1,0)},anchor=south east}, compat=1.3,
  xmin=0, xmax=4,ymin=0,ymax=0.5,
  xlabel= {Time $T$},
  ylabel= {$\mathbb{E} ( Y_T^2 )$},
  tick label style={/pgf/number format/fixed},
  legend cell align=left]
 \addplot[thick,dashed,color=black,mark=none] table [x index=0, y index=1]{./sL_mkenergy_pde100_dt1.0e-03_data.txt};
 \addplot[thick,dotted,color=black,mark=none] table [x index=0, y index=1]{./sL_mkenergy_pde200_dt1.0e-03_data.txt};
 \addplot[thick,dashdotted,color=black,mark=none] table [x index=0, y index=1]{./sL_mkenergy_pde300_dt1.0e-03_data.txt};
 \addplot[thick,color=black,mark=none] table [x index=0, y index=1]{./sL_mkenergy_pde350_dt1.0e-03_data.txt};
 \legend{$L=0.75$,$L=1.5$,$L=2.25$,$L=2.625$}
\end{axis}
\node at (0.5cm,5cm) {\textcolor{black}{(b)}};
\end{tikzpicture}
 \begin{tikzpicture}[scale=0.8125]
\begin{axis}[legend style={at={(1,0)},anchor=south east}, compat=1.3,
  xmin=0, xmax=4,ymin=0,ymax=0.6,
  xlabel= {Time $T$},
  ylabel= {$\mathbb{P} \left ( |Z_T| = P_Y \right )$},
  tick label style={/pgf/number format/fixed},
    legend cell align=left]
\addplot[thick,solid,color=black!20!white, mark=none] table [x index=0, y index=1]{./plasticstate_pde25_dt1.0e-04_data.txt};
\addplot[thick,solid,color=black!35!white, mark=none] table [x index=0, y index=1]{./plasticstate_pde50_dt1.0e-04_data.txt};
\addplot[thick,solid,color=black!50!white, mark=none] table [x index=0, y index=1]{./plasticstate_pde100_dt1.0e-04_data.txt};
\addplot[thick,solid,color=black!65!white, mark=none] table [x index=0, y index=1]{./plasticstate_pde200_dt1.0e-04_data.txt};
\addplot[thick,solid,color=black!80!white, mark=none] table [x index=0, y index=1]{./plasticstate_pde400_dt1.0e-04_data.txt};
\addplot[thick,solid,color=black!100!white,mark=none] table [x index=0, y index=1]{./plasticstate_pde800_dt1.0e-04_data.txt};
    \legend{$\tilde I = \tilde J = 25$,$\tilde I = \tilde J = 50$,$\tilde I = \tilde J = 100$,$\tilde I = \tilde J = 200$,$\tilde I = \tilde J = 400$,$\tilde I = \tilde J = 800$}
\end{axis}
\node at (0.5cm,5cm) {\textcolor{black}{(c)}};
\end{tikzpicture}
\begin{tikzpicture}[scale=0.8125]
\begin{axis}[legend style={at={(1,0)},anchor=south east}, compat=1.3,
  xmin=0, xmax=4,ymin=0,ymax=0.5,
  xlabel= {Time $T$},
  ylabel= {$\mathbb{E} ( Y_T^2 )$},
  tick label style={/pgf/number format/fixed},
  legend cell align=left]
\addplot[thick,solid,color=black!20!white, mark=none] table [x index=0, y index=1]{./mkenergy_pde25_dt1.0e-04_data.txt};
\addplot[thick,solid,color=black!35!white, mark=none] table [x index=0, y index=1]{./mkenergy_pde50_dt1.0e-04_data.txt};
\addplot[thick,solid,color=black!50!white, mark=none] table [x index=0, y index=1]{./mkenergy_pde100_dt1.0e-04_data.txt};
\addplot[thick,solid,color=black!65!white, mark=none] table [x index=0, y index=1]{./mkenergy_pde200_dt1.0e-04_data.txt};
\addplot[thick,solid,color=black!80!white, mark=none] table [x index=0, y index=1]{./mkenergy_pde400_dt1.0e-04_data.txt};
\addplot[thick,solid,color=black!100!white, mark=none] table [x index=0, y index=1]{./mkenergy_pde800_dt1.0e-04_data.txt};
    \legend{$\tilde I = \tilde J = 25$,$\tilde I = \tilde J = 50$,$\tilde I = \tilde J = 100$,$\tilde I = \tilde J = 200$,$\tilde I = \tilde J = 400$,$\tilde I = \tilde J = 800$}
\end{axis}
\node at (0.5cm,5cm) {\textcolor{black}{(d)}};
\end{tikzpicture}
\caption{Empirical convergence study for $\mathbb{P} (
  |Z_T| = P_Y )$ and $\mathbb{E} ( Y_T^2 )$, $T \in [0,4]$. In all
  cases, $c_0 = 1, k=1, P_Y = 0.25$.  PDE solutions of the probability
  of plastic state and the mean kinetic energy are shown for different
  values of the truncation bound $L$ (while $\delta y = 0.0075$ is
  kept constant) in the $y$-direction in (a) and (b), and different
  values of $\tilde{I} = \tilde{J}$ in (c) and (d) while keeping $L = 3.0$. 
  \label{fig:conv_study}}
\end{figure}

We first present an empirical study on the convergence of the
probability of the plastic state \eqref{epp_plasticstate} and the mean
kinetic energy \eqref{epp_mkenergy}. The results presented in Figure
\ref{fig:conv_study} provide insight into the dependence of the PDE
solution on the domain truncation $L$ and on the number of
discretization points. Note in particular, that the convergence of
$\mathbb{E} ( Y_T^2 )$ requires a sufficiently large value of $L$,
which is due to the fact that this quantity involves squared $y$-values. Moreover, we also observe that a
sufficiently fine mesh with $\tilde I=\tilde J\ge 100$ is required to
properly resolve the PDE problem.\\

\begin{figure}[h!]
\centering
\begin{tikzpicture}[scale=0.8125]
\begin{axis}[legend style={at={(1,0)},anchor=south east}, compat=1.3,
  xmin=0, xmax=4,ymin=0,ymax=0.5,
  xlabel= {Time $T$},
  ylabel= {$\mathbb{P} \left ( |Z_T| = P_Y \right )$.}]
\addplot[thick,dashed,color=black,mark=none,mark size=1pt] table [x index=0,
  y index=1]{./plasticstate_mc6dt5_data.txt};
\addlegendentry{\#MC = $10^6$, $\delta t = 10^{-5}$}
\addplot[thick,color=black, mark=none] table [x index=0, y index=1]{./plasticstate_pde800_dt1.0e-04_data.txt};
\addlegendentry{$\tilde I = \tilde J = 800$, $\delta t = 10^{-4}$}
\end{axis}
\node at (0.5cm,5cm) {\textcolor{black}{(a)}};
\end{tikzpicture}\hspace*{.4cm}
\begin{tikzpicture}[scale=0.8125]
\begin{axis}[legend style={at={(1,0)},anchor=south east}, compat=1.3,
  xmin=0, xmax=4,ymin=0,ymax=0.4,
  xlabel= {Time $T$},
  ylabel= {$\mathbb{P} \left ( |Z_T| = P_Y, |Z_{T+h}| = P_Y \right )$}]
\addplot[thick,dashed,color=black,mark=none,mark size=1pt] table [x index=0, y index=1]{./corr_plasticstate_mc6dt5_data.txt}; 
\addlegendentry{\#MC = $10^6$, $\delta t = 10^{-5}$}
\addplot[thick,color=black,mark=none,mark size=1pt] table [x index=0, y index=1]{./corr_plasticstate_pde800_dt1.0e-04_data.txt};
\addlegendentry{$\tilde I = \tilde J = 800$, $\delta t = 10^{-4}$}
\end{axis}
\node at (0.5cm,5cm) {\textcolor{black}{(b)}};
\end{tikzpicture}
\begin{tikzpicture}[scale=0.8125]
\begin{axis}[legend style={at={(1,0)},anchor=south east}, compat=1.3,
  xmin=0, xmax=4,ymin=0,ymax=0.5,
  xlabel= {Time $T$},
  ylabel= {$\mathbb{E} Y_T^2$}]
\addplot[thick,dashed,color=black,mark=none,mark size=1pt] table [x index=0, y index=1]{./mkenergy_mc6dt5_data.txt}; 
\addlegendentry{\#MC = $10^6$, $\delta t = 10^{-5}$}
\addplot[thick,color=black, mark=none] table [x index=0, y index=1]{./mkenergy_pde800_dt1.0e-04_data.txt};
\addlegendentry{$\tilde I = \tilde J = 800$, $\delta t = 10^{-4}$}
\end{axis}
\node at (0.5cm,5cm) {\textcolor{black}{(c)}};
\end{tikzpicture}
\hspace*{.4cm}
\begin{tikzpicture}[scale=0.8125]
\begin{axis}[legend style={at={(1,0)},anchor=south east}, compat=1.3,
  xmin=0, xmax=4,ymin=0,ymax=0.5,
  xlabel= {Time $T$},
  ylabel= {$\mathbb{E} Y_T^2 Y_{T+h}^2$}]
\addplot[thick, dashed,color=black,mark=none,mark size=2pt] table [x index=0, y index=1]{./corr_mkenergy_mc6dt5_data.txt}; 
\addlegendentry{\#MC = $10^6$, $\delta t = 10^{-5}$}
\addplot[thick,color=black,mark=none,mark size=2pt] table [x index=0, y index=1]{./correlation_mkenergy_pde800_dt1.0e-04_data.txt};
\addlegendentry{$\tilde I = \tilde J = 800$, $\delta t = 10^{-4}$}
\end{axis}
\node at (0.5cm,5cm) {\textcolor{black}{(d)}};
\end{tikzpicture}
\begin{tikzpicture}[scale=0.8125]
\begin{axis}[legend style={at={(0,1)},anchor=north west}, compat=1.3,
  xmin=0, xmax=4,ymin=0,ymax=2,
  xlabel= {Time $T$},
  ylabel= {$\mathbb{E} \left ( \Delta_T - \mathbb{E} \Delta_T \right )^2$}]
\addplot[thick,dashed,color=black,mark=none,mark size=1pt] table [x index=0, y index=1]{./def_total_variance_mc6dt5_data.txt}; 
\addlegendentry{\#MC = $10^6$, $\delta t = 10^{-5}$}
\addplot[thick,color=black,mark=none,mark size=1pt] table [x index=0, y index=1]{./var_total_pde800_dt1.0e-04_data.txt};
\addlegendentry{$\tilde I = \tilde J = 800$, $\delta t = 10^{-4}$}
\end{axis}
\node at (0.5cm,3.8cm) {\textcolor{black}{(e)}};
\end{tikzpicture}
\hspace*{.4cm}
\begin{tikzpicture}[scale=0.8125]
\begin{axis}[legend style={at={(0,1)},anchor=north west}, compat=1.3,
  xmin=0, xmax=4,ymin=0,ymax=2,
  xlabel= {Time $T$},
  ylabel= {$\mathbb{E} \left ( X_T - \mathbb{E} X_T \right )^2$}]
\addplot[thick,dashed,color=black,mark=none,mark size=1pt] table [x index=0, y index=1]{./def_plastic_variance_mc6dt5_data.txt}; 
\addlegendentry{\#MC = $10^6$, $\delta t = 10^{-5}$}
\addplot[thick,color=black,mark=none,mark size=1pt] table [x index=0, y index=1]{./def_plastic_variance_pde1000dt3_data.txt};
\addlegendentry{$\tilde I = \tilde J = 800$, $\delta t = 10^{-4}$}
\end{axis}
\node at (0.5cm,3.8cm) {\textcolor{black}{(f)}};
\end{tikzpicture}
\caption{Comparison between the PDE (solid) and MC (dashed) solutions
  for $T \in [0,4]$ for the elasto-plastic problem.  (a): probability of the plastic state
  $\mathbb{P} \left ( |Z_T| = P_Y \right )$.  (b): a correlation
  structure of the probability of the plastic state $\mathbb{P} \left
  ( |Z_T| = P_Y, |Z_{T+h}| = P_Y \right ),h= 0.2$.  (c): mean
  kinetic energy $\mathbb{E} Y_T^2$.  (d): a correlation
  structure of the mean kinetic energy $\mathbb{E} Y_T^2 Y_{T+h}^2$.
  (e): the variance of the plastic deformation $\mathbb{E}
  \left ( \Delta_T - \mathbb{E} \Delta_T \right )^2$.  (f):
  the variance of the total deformation $\mathbb{E} \left ( X_T -
  \mathbb{E} X_T \right )^2$.}
\label{fig:epp_results}
\end{figure}
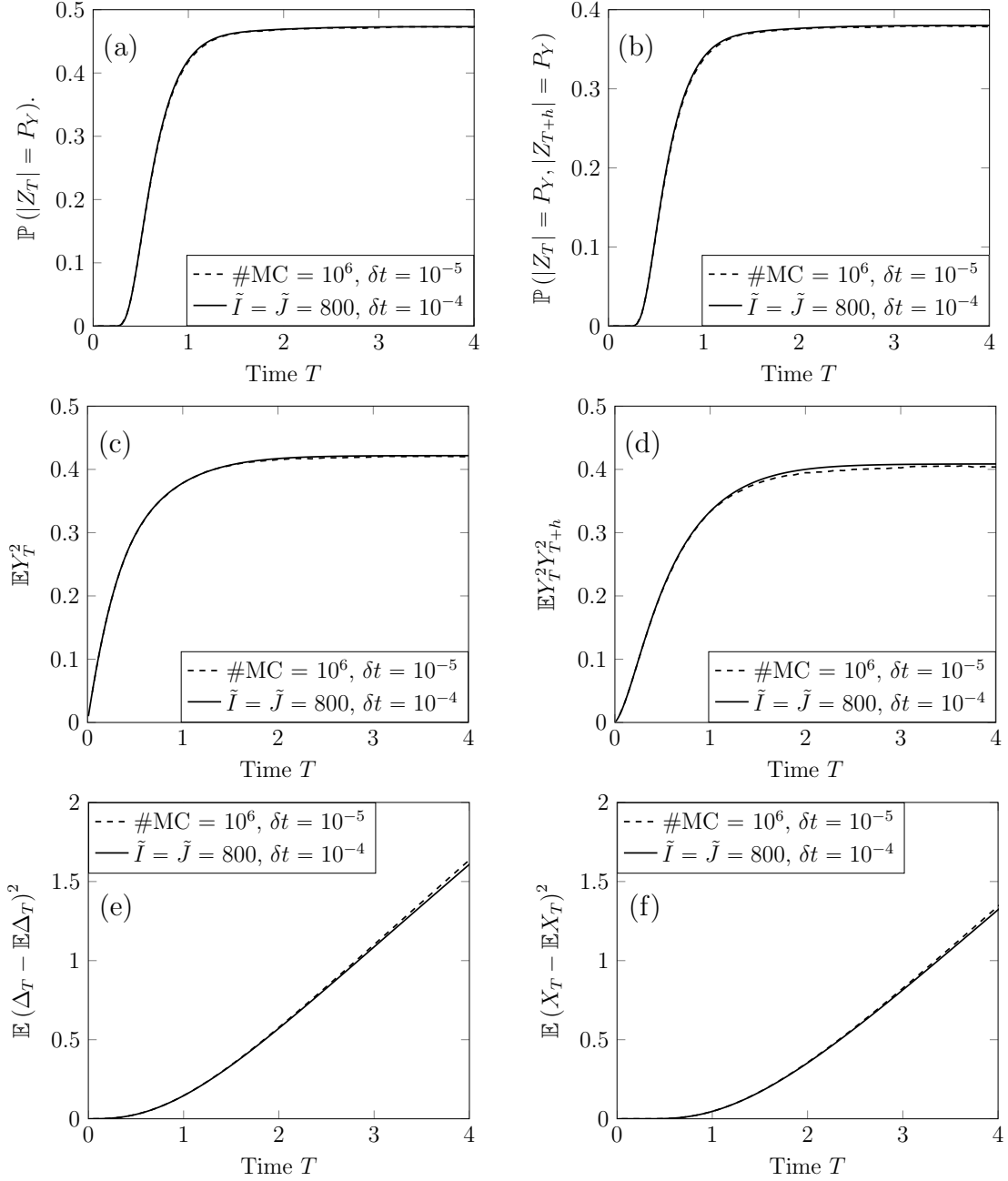

\noindent
Next, in Figure~\ref{fig:epp_results}, we present systematic
numerical comparisons between the PDE and the probablistic Monte Carlo (MC)
approach. In particular, we show comparisons for the quantities
\eqref{epp_plasticstate}, \eqref{epp_mkenergy}, \eqref{epp_variance}, \eqref{epp_correlation} \lmxxx{and \eqref{epp_crossings}}. We have chosen discretization
parameters, for which we found small approximation errors in our
previous tests summarized in Figure~\ref{fig:conv_study}.  As can be
seen, the results found from the PDE solution closely track the
results from the probabilistic simulation. However, some quantities are harder to approximate than others and it appears that for the correlation structure in the kinetic energy, as shown in Figure~\ref{fig:epp_results}(d), there is a slight discrepancy between the PDE and probabilistic results. This may be related to the fact that the quantity is quartic in the velocity.
From Figure~\ref{fig:epp_results}, we observe as expected from theory
(see \cite{BT08}) that $(Z_T,Y_T)$ has a unique invariant probability
measure. Thus, when $T$ becomes large, for any well-behaved function
$f$, $\mathbb{E} f(Z_T,Y_T)$ becomes constant and $\mathbb{E} \left (
\int_0^T f(Z_t,Y_t) \d t \right )^2$ obeys linear growth with respect
to $T$. To the best of our knowledge, except quantities of type $A$
for large time \cite{BMPT09,BFMY15}, quantities of type $A,B,C$ or
$A',B',C'$ (both transient and stationary case) have not been computed
using PDE formulations in such a general framework previously.

\subsection{Numerical results for the obstacle problem}
Next, we present numerical results for the stochastic obstacle
problem. As above, for the PDE formulation \eqref{uop}, \eqref{vop}
and \eqref{wop} we use a finite-difference scheme in space and the
implicit Euler method in time. The main difference compared to the
discretization for the elasto-plastic problem presented in
Section~\ref{ssec:num1} is in the discretization of the non-standard
boundary conditions to be satisfied at the
red points in Figure~\ref{fig:grid}.
We recall that the boundary conditions on $D_T^\pm$ in the PDEs for the elasto-plastic and obstacle problems 
are different because they reflect different boundary behaviors for the underlying stochastic processes.
In the elasto-plastic problem, the boundary condition (plastic phases) solves 
a boundary value problem whose boundary data is also a part of the problem.
In the obstacle problem, the boundary condition (impacts with the obstacle) consists in identifying the values of the solution on $D_T^\pm$ to those on $x = \pm 1, \mp y >0$.    
\begin{align*}
& u_{1,j} = c_j u_{1, j_e} + (1-c_j) u_{1, j_e+1}  \quad \mbox{and} \quad v_{1,j} = c_j v_{1, j_e} + (1-c_j) v_{1, j_e+1}  \quad\text{for } 1 \le j\le \tilde J,\\
& u_{I,j} = c_j u_{I, j_e} + (1-c_j) u_{I, j_e+1}  \quad \mbox{and} \quad v_{I,j} = c_j v_{I, j_e} + (1-c_j) v_{I, j_e+1} \quad \text{for } \tilde J+2\le j\le J,
\end{align*}
where
\[
j_e \triangleq 1 + \left [ \frac{1}{\delta y} (L_y - e y_j) \right ]\quad \mbox{and} 
\quad c_j \triangleq \frac{y_{j_e+1}+e y_j}{\delta y}.
\]
We note that $y_j$ is defined to be on the grid, but, in general, the quantity $e y_j$ will not correspond to a grid point for $e \in (0,1)$.
Consequently, the value imposed in $u_ {1, j}$ (or in $u_ {I, j}$) is a value interpolated between the values of $u$ at the two nearest neighbors of $-ey_j$ in the grid. Here we use the notation $[ a ]$ for the integer part of $a$. In the purely elastic case $e = 1$, $j_1 = J - (j-1)$ and
$c_j = 1$ whereas in the purely inelastic case $e = 0$, $j_0
\equiv \tilde{J} + 1$ and $c_j = 1$.\\

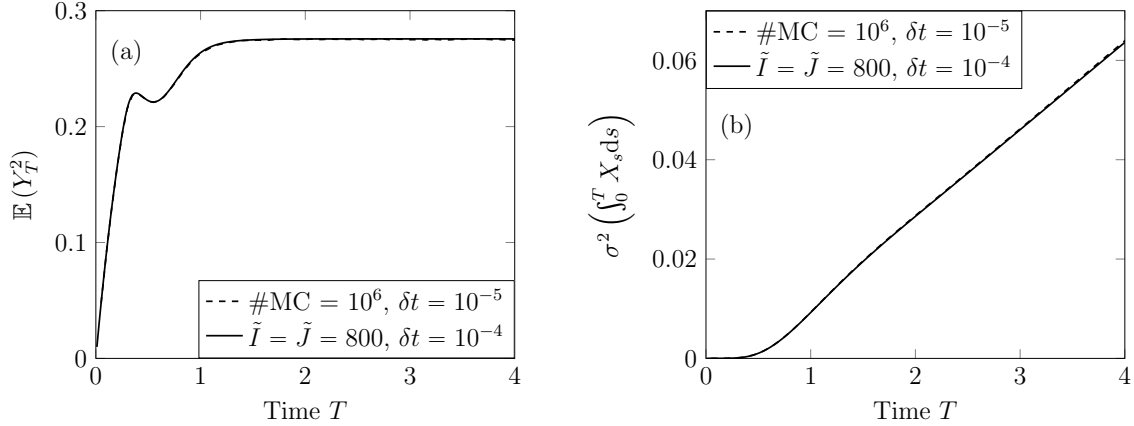
\begin{figure}[h]  
\centering
\begin{tikzpicture}[scale=0.8125]
\begin{axis}[legend style={at={(1,0)},anchor=south east}, compat=1.3,
  xmin=0, xmax=4,ymin=0,ymax=0.3,
  xlabel= {Time $T$},
  ylabel= {$\mathbb{E} \left ( Y_T^2 \right )$}]
\addplot[thick,dashed,color=black,mark=none,mark size=1pt] table [x index=0, y index=1]{./obstaclemkenergy_dt1.0E-04_MC_1.0E+06_data.txt}; 
\addlegendentry{\#MC = $10^6$, $\delta t = 10^{-5}$}
\addplot[thick,color=black,mark=none,mark size=1pt] table [x index=0, y index=1]{./mkenergyobstacle_pde800_dt1.0e-04_data.txt};
\addlegendentry{$\tilde I = \tilde J = 800$, $\delta t = 10^{-4}$}
\node at (0.5cm,5cm) {\textcolor{black}{(a)}};
\end{axis}
\end{tikzpicture}
\hspace*{.4cm}
\begin{tikzpicture}[scale=0.8125]
\begin{axis}[legend style={at={(0,1)},anchor=north west}, compat=1.3,
  xmin=0, xmax=4,ymin=0,ymax=0.07,
  xlabel= {Time $T$},
  ylabel= {$\sigma^2 \left ( \int_0^T X_s \d s \right )$},
  scaled ticks=false, tick label style={/pgf/number format/fixed}]
\addplot[thick,dashed,color=black,mark=none,mark size=1pt] table [x index=0, y index=1]{./var_obstacle_def_dt1.0E-04_MC_1.0E+06_data.txt}; 
\addlegendentry{\#MC = $10^6$, $\delta t = 10^{-5}$}
\addplot[thick,color=black,mark=none,mark size=1pt] table [x index=0, y index=1]{./var_obstacle_pde500_dt1.0e-03_data.txt};
\addlegendentry{$\tilde I = \tilde J = 800$, $\delta t = 10^{-4}$}
\node at (0.5cm,3.8cm) {\textcolor{black}{(b)}};
\end{axis}
\end{tikzpicture}
\caption{Obstacle problem with $e = 0.5$. Comparison between the PDE (solid) and MC (dashed) solutions for $T \in [0,4]$. 
 (a): mean kinetic energy
  $\mathbb{E} \left ( Y_T^2 \right )$.  (b): the variance of the
  integral of the total deformation $\sigma^2 \left ( \int_0^T X_s \d
  s \right )$.}
\label{fig:op_results}
\end{figure}

\noindent
In Figure \ref{fig:op_results}, we present numerical comparisons
between the PDE and probabilistic MC approaches for
\eqref{op_mkenergy} and \eqref{op_variance}. As can be seen, the
solution of the PDE approach agrees well with the solution of the
probabilistic simulation. We also observe that when $T$ becomes large,
the expectation (left plot) becomes constant and the variance (right
plot) grows linearly. While we expect that ergodicity holds for the
obstacle problem, to the best of our knowledge this has not been proven
in the literature. Additionally, as far as we are aware of, quantities
of type $A,B,C$ or $A',B',C'$
(both transient and stationary case) have not been computed using
PDEs in such a general framework before. Observe that \eqref{op_mkenergy} and \eqref{op_variance} are of type $A$ and $A'$ when $T$ is finite.

\subsection{Long-time behavior}
\noindent
Next, we present results for the solution of the stationary versions
of \eqref{uep} combined with \eqref{wep}, and of \eqref{uop} combined
with \eqref{wop} where $\lambda$ is chosen small enough. We remind the
reader that $\lambda u$ approximates a quantity of the form $\lim_{ t
  \to \infty} \mathbb{E} f(Z_t, Y_t)$ when $\lambda$ is small. We use
a similar finite-difference procedure as above for the time-dependent
case to obtain a linear system of the form
\begin{equation}
\label{algorithm_stationary}
(\lambda I + M) \boldsymbol{u} = \boldsymbol{f}.
\end{equation}
\noindent Again, we first compare results computed using the
discretized PDE approach \eqref{algorithm_stationary} with results
based on Monte Carlo simulations.  For the elasto-plastic PDE
discretization, we use $\tilde I = \tilde J = 800$, i.e., the overall
number of spatial unknowns is $(2\tilde I+1)(2\tilde J+1) =
2,563,201$, and we choose $\lambda=10^{-3}$.  For the Monte Carlo
simulation, we use $\delta t = 10^{-5}$. We present the comparison
between the PDE and the Monte Carlo (MC) approach in Table
\ref{tab1}. For the obstacle problem, the spatial discretization uses
$\tilde I = \tilde J = 500$, amounting to overall $1,002,001$
unknowns. We again use $\lambda=10^{-3}$. Results of the comparison
between the PDE and the Monte Carlo approach are given in Table
\ref{tab2}. As can be seen, for both the obstacle as well as the
elasto-plastic problem, the results obtained with Monte Carlo are
close to the PDE results. To study the interplay between the time discretization and the Monte
Carlo errors, we considered the case of the growth rate of the
variance related to the plastic deformation in the elasto-plastic
problem. \lm{Approximations of \eqref{epp_plasticstate} and \eqref{epp_variance} are shown in Table
\ref{tabemper} for $P_Y = 0.25$ and for different values of $\delta t$
and $T$. The approximation of \eqref{epp_plasticstate} is given by \eqref{approxproba} with $g= 0$ and $f(y,z) = \mathbf{1}_{ \{ |z| = P_Y\} }$ and the approximation of \eqref{epp_variance} is given by
$$\frac{1}{T M} \sum_{m=1}^M \left (
\sum_{n=0}^{N_{\delta t}-1} g(Z_n^m,Y_n^m) (t_{n+1}-t_n) \right )^2,
$$
where $g(y,z) = y \mathbf{1}_{ \{ |z| = 1\} }$. In both approximations, $(Z_n^m,Y_n^m)$
satisfies the probabilistic numerical scheme for \eqref{svi_ep} and the Monte Carlo sample size is $10^5$. 
For the approximation of \eqref{epp_plasticstate} we observe a relatively fast convergence towards a constant with respect to $T \geq 3$ whereas for the approximation of \eqref{epp_variance} a sufficiently large value of $T \geq 150$ is required to see the convergence.}
For both cases, a sufficiently small value of $\delta t = 10^{-4}$ is required. Empirically, the data indicates that the error is roughly halved as the time step is divided by 10.
\begin{table}[h]
\centering
\caption{Results using the MC approach with different $T$ and $\delta t$ 
for the approximation of the probability \eqref{epp_plasticstate} (left) and the growth rate of the plastic deformation \eqref{epp_variance} (right). Here $P_Y = 0.25$. 
\label{tabemper}}
\begin{tabular}{|c|c|c|c|}
$T$ & $\delta t = 10^{-2}$ & $ \delta t = 10^{-3}$ & $\delta t  = 10^{-4}$ \\ \hline
2   & 0.518  &   0.487   &   0.474\\ 
4  & 0.519  &   0.492  &    0.478\\
8  & 0.519  &   0.489  &    0.477\\
16  & 0.520  &  0.487 &   0.481\\
32  & 0.518  &  0.491  &   0.478\\
64  & 0.519  &  0.490 &   0.476\\
\end{tabular}
\hspace{0.5cm}
\begin{tabular}{|c|c|c|c|}
$T$ & $\delta t = 10^{-2}$ & $ \delta t = 10^{-3}$ & $\delta t  = 10^{-4}$ \\ \hline
10  & 0.489  &    0.464  &    0.456\\
20  & 0.532  &    0.505  &    0.493\\
40  & 0.558  &    0.528   &   0.510\\
80  & 0.567  &    0.538    &  0.521\\
160  & 0.575  &   0.548    &  0.524\\
200  & 0.574  &   0.547    &  0.526
\end{tabular}
\end{table}
\begin{table}[h]
\centering
\caption{Results using PDE and MC method approaches for the elasto-plastic
  problem \eqref{svi_ep} for different plastic yields $P_Y$:
  Probability \eqref{epp_plasticstate} of plastic state for
  $T\to\infty, (T \geq 5)$. Asymptotic growth rate \eqref{epp_variance} of the
  plastic/total deformation for $T\to\infty, (T \geq 150)$.\label{tab1}}
\begin{tabular}{c|l|llllllllll}
&$P_Y$ & 0.1 & 0.2 & 0.3 & 0.4 & 0.5 & 0.6 & 0.7 & 0.8 & 0.9 & 1.0 \\ \hline
\eqref{epp_plasticstate}&\textbf{PDE} &  0.638 & 0.521 & 0.430 & 0.354 &  0.289 & 0.234  & 0.187 & 0.148  & 0.115  & 0.088 \\
$T\!\to\!\infty$&\textbf{MC} &  0.639 & 0.521 & 0.429 & 0.352 & 0.286 & 0.230 & 0.184 &
0.144 & 0.112 &   0.085\\ \hline
\eqref{epp_variance}&\textbf{PDE}  &  0.776  & 0.589 & 0.440 & 0.325 & 0.238  & 0.173  & 0.125 & 0.089 & 0.064  & 0.045 \\
$T\!\to\!\infty$&\textbf{MC} &  0.799 & 0.607 & 0.457  & 0.337 & 0.248  & 0.182  & 0.132 & 0.096  &  0.070 &  0.050 \\ 
\end{tabular}
\end{table}
\begin{table}[h]
\centering
\caption{Results using PDE and MC method approaches for the obstacle
  problem \eqref{svi_obstacle} for different position of the obstacle $P_O$:
  Mean kinetic energy \eqref{op_mkenergy} for $T\to\infty, (T \geq 5)$. 
  Asymptotic growth rate \eqref{epp_variance} of the
  integral of the deformation for $T\to\infty, (T \geq 150)$.\label{tab2}}
\begin{tabular}{c|l|llllllllll}
&$P_Y$ & 0.1 & 0.2 & 0.3 & 0.4 & 0.5 & 0.6 & 0.7 & 0.8 & 0.9 & 1.0 \\ \hline
\eqref{op_mkenergy} & \textbf{PDE} &  0.179 & 0.250 & 0.298 & 0.347 &  0.364 & 0.388  & 0.409 & 0.426  & 0.441  & 0.453\\
$T\!\to\!\infty$ & \textbf{MC}               &  0.179 & 0.250 & 0.297 & 0.334 & 0.364 & 0.389 & 0.408 & 0.423 & 0.437 &   0.450\\ \hline
\eqref{op_variance}&\textbf{PDE}  &  0.00148  & 0.0097  & 0.028  & 0.062 & 0.103  & 0.158  & 0.223 & 0.294 & 0.370  & 0.447 \\
$T\!\to\!\infty$&\textbf{MC} & 0.00149 & 0.0097  & 0.028 & 0.060 & 0.103 & 0.158  & 0.222  & 0.295 & 0.370  &  0.447  \\ 
\end{tabular}
\end{table}

\section{Comparison with \cite{BMPT09,BMY12,BFMY15} for a white noise EPPO at large time}
\noindent In this section, we compare the approach proposed here with 
previous techniques employed for a white noise EPPO.
\subsection{Approach in \cite{BMPT09} for the elasto-plastic problem}
If the values of the function $u$ were known at $z=P_Y, y=0^+$ and $z=-P_Y, y=0^-$, 
the stationary version of \eqref{uep} would be a degenerate elliptic PDE with
a standard Dirichlet boundary condition. 
This PDE could be solved by a standard numerical approach. 
However, the challenge in solving this problem
resides in the fact that these values are not input data of the
problem but part of the solution. This makes it a non-standard problem and more challenging to solve. As a remedy, a superposition approach has
been proposed in \cite{BMPT09}. To ensure continuity of $u$ at the
points $(P_Y,0)$ and $(-P_Y,0)$, the linearity of $\eqref{uep}$ allows to compute the solution as a linear
combination of the following three local problems:
\[
\lambda v + A_2 v  =  g \hbox{ in  } \: D_2, 
\quad \lambda v + B_{2,+} v  =  g_{+} \hbox{ in  }\: D_2^+,
\quad \lambda v + B_{2,-} v  =  g_{-}\hbox{ in  }\: D_2^-,
\]
\begin{center}
with $v(P_Y,0)=0, v(-P_Y,0)=0$,
\end{center}
\[
\lambda \pi^+ + A_2 \pi^+  =  0 \hbox{ in  }   \: D_2,
\quad \lambda \pi^+ + B_{2,+} \pi^+ = 0 \hbox{ in  } \: D_2^+,
\quad \lambda \pi^+ + B_{2,-} \pi^+ = 0 \hbox{ in  }   \: D_2^-,
\]
\begin{center}
with $\pi^+(P_Y,0)=1, \pi^+(-P_Y,0)=0$,
\end{center}
\[
\lambda \pi^- + A_2 \pi^-  =  0 \hbox{ in  }   \: D_2,
\quad \lambda \pi^- + B_{2,+} \pi^- = 0 \hbox{ in  } \: D_2^+,
\quad \lambda \pi^- + B_{2,-} \pi^- = 0 \hbox{ in  }   \: D_2^-,
\]
\begin{center}
with $\pi^-(P_Y,0)=0, \pi^-(-P_Y,0)=1$.
\end{center}
Then, finding a continuous solution $u$ amounts to finding scalar values $u_+$ and
$u_-$ such that $u = v + u_+ \pi^+ + u_- \pi^-$ is continuous in
$(-P_Y,0)$ and $(P_Y,0)$, which requires to solve the following linear
$2\times 2$ system:
\[
\Pi
\begin{pmatrix}
u_+ \\
u_-
\end{pmatrix}
\!\!=\!\!
\begin{pmatrix}
v(P_Y,0^-)-v(P_Y,0^+)\\
v(-P_Y,0^-)-v(-P_Y,0^+)
\end{pmatrix},
\]
where
\[
\Pi \triangleq 
\begin{pmatrix}
\pi^+(P_Y,0^+) - \pi^+(P_Y,0^-) & \pi^-(P_Y,0^+) - \pi^-(P_Y,0^-)\\
\pi^+(-P_Y,0^+) - \pi^+(-P_Y,0^-) & \pi^-(-P_Y,0^+) - \pi^-(-P_Y,0^-)
\end{pmatrix}.
\]
\begin{remark}
It is important to note that this 
technique in \cite{BMPT09} is based on the
superposition of local PDEs and is thus specific to the case of a SVI
modelling an elasto-perfectly-plastic oscillator (EPPO) excited by
white noise. It has no natural extension to more general EP or obstacle problems.
Indeed, as explained above, solving the KE with non-standard boundary
conditions corresponding to an EPPO excited by white noise
reduces to solving two linear equations for two unknown scalars. 
For more general problems as targeted here, one must find unknown functions, and thus the
superposition method cannot be employed straightforwardly.
This superposition technique has a probabilistic interpretation in terms
of novel notion of short cycles as defined in \cite{BM12}.
\end{remark}

\subsection{Approach in \cite{BFMY15,BMY12} for the growth rate of the variance of $\Delta_t$ at large time}
\noindent In \cite{BMY12}, the growth rate of the variance of $\Delta_t$ at large time has been characterized. 
Relying on \eqref{svi_ep}, the authors proposed a novel and simple formulation of the evolution of the system in terms of stopping times in order to identify a repeating pattern in the trajectory, namely long cycles. This concept is summarized next.
\begin{definition}
A long cycle of the solution $(Z_t,Y_t)$ of \eqref{svi_ep} is a path, enclosed by the stopping times $\tau_0$ and $\tau_1$ defined below, starting and ending in one of the two points $\{ (-P_Y,0), (P_Y,0)\}$ which has touched the other point at least once. Similarly, a half long cycle is a path enclosed by the stopping times $\tau_0$ and $s_0$, or by $s_0$ and $\tau_1$. In a recursive way, a sequence of stopping times $\tau_n$ can be defined where $\tau_n$ is the time when the $n$-th long cycle ends.
\end{definition}
\noindent With the notation 
\[
\tau_0 \triangleq  \inf \{ t > 0, \quad Y_t = 0 \quad \mbox{and} \quad \vert Z_t \vert = P_Y \},
\] 
$\delta\triangleq \textup{sign}(Z_{\tau_0})$, which labels the first boundary reached by the process $(Z_t,Y_t)$, and 
\[
\begin{cases}
& s_0 \triangleq  \inf \{ t > \tau_0, \quad Y_t = 0 \quad \mbox{and} \quad Z_t = -\delta P_Y \},\\
& \tau_1 \triangleq  \inf \{ t > s_0, \quad Y_t = 0 \quad \mbox{and} \quad Z_t = \delta P_Y \},
\end{cases}
\]
they obtained a probabilistic expression for the coefficient of the growth rate of the variance of $\Delta_t$ as follows: 
\begin{equation}
\label{FLC}
\lim_{t \to \infty} \frac{ \sigma^2(\Delta_t)}{t} = \mu \gamma^2,
\end{equation}
where
\begin{equation}\label{eq:DefMuGamma}
\gamma^2 \triangleq \mathbb{E} \left ( \Delta_{\tau_1}-\Delta_{\tau_0} \right)^2 \quad \mbox{and} \quad \mu \triangleq \frac{1}{\mathbb{E} (\tau_1 - \tau_0)}.
\end{equation}
In \cite{BFMY15}, using a Fourier transform approach, an analytic framework for $\mu$ and $\gamma^2$ has been proposed. 

\begin{remark}
It is important to note that the technique in \cite{BMY12,BFMY15} of splitting the trajectory in terms of long cycle (an identically independent distributed repeating pattern) 
is specific to the case of a SVI modelling an EPPO excited by white noise. It has no natural extension to more general EP, obstacle or colored noise problems.
Moreover, it cannot be employed for short durations.
\end{remark}

\section{Conclusions and perspectives}
\noindent 
\lmxxx{PDEs with non-standard boundary conditions for FKf} that describe non-smooth stochastic processes have been (formally) derived and numerically solved. Important examples from engineering have been successfully treated. The present work will pave the way for promising new research. Indeed, the technique presented in this paper is straightforwardly generalizable to a broad range of cases where $F$ and $G$ are time-dependent and to problems driven by colored noise. This framework can also be extended to study Power Spectral Densities \textcolor{black}{by using a system of PDEs with non-standard boundary conditions. 
Indeed, considering for instance the impact problem, if a complex valued function $u \in C^\star (\mathbb{R}^2 \times [0,T])$ and a real valued function $v \in C^\star (\mathbb{R}^2 \times [0,T])$ 
satisfy
\begin{align*}
\label{cu}
& \frac{\partial u}{\partial t} + L u - \mathbf{i} \xi u = -g(x,y), \quad u(T) = 0 
& \frac{\partial v}{\partial t} +  L v  = - \left | \nabla u \right |^2, \quad v(T) = 0
\end{align*}
together with the boundary conditions for the impact problem,
then we have 
\[
\mathbb{E} \left ( \left | \int_0^T  g(X_t,Y_t) e^{-\mathbf{i} \omega t} \textup{d} t  \right |^2 \right  ) = v(x,y,0) + |u(x,y,0)|^2.
\]
Here, $L$ is the infinitesimal generator of $(X,Y)$ away from the boundary.}

\appendix
\section{Probabilistic simulation}
\label{algo_svi_obstacle}
\noindent Below a detailed implementation of the probabilistic simulation for \eqref{svi_obstacle}.
\noindent For $T>0, N \in \mathbb{N}$ and $\delta t \triangleq \frac{T}{N}$.
We set $\Sigma \in \mathbb{R}^{2 \times 2}$ such that 
\begin{equation}
\Sigma \Sigma^{T}=
\begin{pmatrix} 
      \sigma_x^2(\delta t) & \sigma_{x,y}(\delta t)\\
      \sigma_{x,y}(\delta t) & \sigma_{y}^2(\delta t)\\
   \end{pmatrix},
\end{equation}
where
\begin{align*}
\sigma_x^2(t) &  =  \dfrac{1}{\omega^2} \int_0^t e^{-c_0s}\sin^2{(\omega s)}ds,\\
\sigma_{y}^2(t) & = \int_0^t e^{-c_0s}\cos^2{(\omega s)}ds - \frac{4 c_0^2}{\omega^2} \int_0^t e^{-c_0s}\sin^2{(\omega s)}ds -  \frac{ c_0}{2 \omega^2}e^{-c_0t}\sin^2{(\omega t)},\\
\sigma_{xy}(t) & =  \frac{1}{2 \omega} \int_0^t e^{-c_0s}\sin{(2 \omega s)}ds -  \frac{c_0}{4 \omega^2}\int_0^t e^{-c_0s}\sin^2{(\omega s )}ds.\\
\end{align*}      
Here, $\omega  \triangleq \dfrac{\sqrt{ 4k - c_0^2}}{2}$, it is assumed that $4k >c_0^2$. 
For every $(x,y) \in \mathbb{R}^2$, we define 
\[
M( \delta t, x,y) \triangleq
\begin{pmatrix} 
m_1(\delta t, x,y)\\
m_2(\delta t, x,y)\\
\end{pmatrix}
\]
where
\begin{eqnarray*}
m_1(t,x,y) & = & e^{\frac{-c_0t}{2}} \lbrace x \cos{(\omega t )} + \frac{1}{\omega}(y+\dfrac{c_0}{2}x)\sin{(\omega t )}\rbrace,\\
m_2(t,x,y)& = & -\frac{c_0}{2}e_x(t,x,y) + e^{-\frac{c_0 t}{2}} \lbrace -\omega x \sin{(\omega t)} + (y+\frac{c_0}{2}x)\cos{(\omega t)}\rbrace.\\
\end{eqnarray*}

\noindent
We can now detail the probabilistic simulation for \eqref{svi_obstacle},
where the main challenge is the incorporation of the obstacles in the
simulation.  Let $(G_{n,m})_{n=0..N,m=1,2}$ be random independent
Gaussian $\mathcal{N}(0,1)$ variables. To compute the $(n+1)$th time step,
we attempt to perform the explicit step
\begin{equation}\label{eq:prob_step}
\begin{pmatrix}
{X}_{n+1}\\
{Y}_{n+1}
\end{pmatrix}
\triangleq 
M \left ( \delta t, 
\begin{pmatrix}
X_n\\
Y_n
\end{pmatrix} \right )
+ \Sigma 
\begin{pmatrix}
G_{n,1}\\
G_{n,2}
\end{pmatrix}.
\end{equation}
If we find that the $(n+1)$st point does not satisfy the obstacle
conditions, for instance since ${X}_{n+1} > 1$, we adjust the
step length to $\theta_{n+1}\delta t$, with $\theta_{n+1} \triangleq
({1-X_n})/({{X}_{n+1} - X_n})$, and set
\begin{equation*}
\quad t_{n+1} \triangleq t_n + \theta_{n+1} \delta t,\quad
X_{n+1} \triangleq 1, 
\quad Y_{n+1} \triangleq -e \left ( (1-\theta_{n+1}) Y_n +
\theta_{n+1} {Y}_{n+1} \right ).
\end{equation*}
An analogous reduction of the step length is used if the full step
computed in \eqref{eq:prob_step} satisfies ${X}_{n+1} < -1$.

\section{Proofs}
\subsection{Bridge from solution of \eqref{Nonsmooth} to PDEs (when $\mathbf{H}=0$) related to $A,B,C$ and $A',B',C'$.}
\noindent This is classic but we put it in Appendix for convenience of the reader.
\begin{proposition}
\label{prop1}
For any function $\phi \in C^\star (\mathbb{R}^2 \times [0,\infty])$ 
satisfying $(\mathcal{A}(\lambda,\phi)$, 
\[
\mathbb{E} \left ( \int_0^t e^{-2 \lambda \tau} \left | \frac{\partial \phi }{\partial y} \right |^2(X_\tau,Y_\tau, \tau) \d \tau \right ) < \infty, \forall t 
\]
the process
\begin{equation}
\label{martingale1}
M_t^{\lambda, \phi} \triangleq 
e^{-\lambda t} \phi(X_t,Y_t,t)
- \int_0^t e^{-\lambda \tau} 
\left ( \frac{\partial \phi}{\partial t} + L \phi  - \lambda \phi \right )(X_\tau,Y_\tau,\tau) \d \tau
\end{equation}
satisfies 
\begin{equation}
\label{martingale2}
M_t^{\lambda,\phi}  = \phi(x,y,0) + \int_0^t  e^{-\lambda \tau} 
\frac{\partial \phi}{\partial y} (X_\tau,Y_\tau, \tau) \d W_\tau
\end{equation}
where $L$ is the differential operator
\[
L  \triangleq \frac{1}{2} \frac{\partial^2 }{\partial y^2} 
+ F \frac{\partial }{\partial x} + G \frac{\partial }{\partial y}.
\]
Thus $M_t^{\lambda, \phi}$ is a martingale under $\mathcal{F}_t \triangleq \sigma \{ W_s, 0 \leq s \leq t \}$. 
Moreover, for any function $\phi' \in C^\star (\mathbb{R}^2 \times [0,\infty])$ satisfying $(\mathcal{A}(\mu,\phi'))$
\begin{equation}
\label{ItoIsometry}
\mathbb{E} \left ( M_T^{\lambda,\phi} M_{T+h}^{\mu,\phi'} \right ) =
\phi \phi'(x,y,0) 
+ \mathbb{E} \left ( \int_0^{T} 
e^{-(\lambda + \mu ) \tau}
\frac{\partial \phi}{\partial y}
\frac{\partial \phi'}{\partial y}(X_\tau,Y_\tau,\tau) \d \tau \right ), \quad \forall T \geq 0, \quad \forall h \geq 0.
\end{equation}
\end{proposition}
\subsection{Proof of Proposition \ref{prop1}}
\begin{proof}
From the assumption on $\phi$ the stochastic integral $\int_0^t  e^{-\lambda \tau} \dfrac{\partial \phi}{\partial y} (X_\tau,Y_\tau, \tau) \d W_\tau$ is well defined. 
Thus, \eqref{martingale2} is obtained by using Ito's formula,
\[
e^{-\lambda t} \phi(X_t,Y_t,t) - \phi(x,y,0) 
= \int_0^t e^{-\lambda \tau} \left ( \frac{\partial \phi}{\partial t} + L \phi  - \lambda \phi \right )(X_\tau,Y_\tau,\tau) \d \tau
+ \int_0^t  e^{-\lambda \tau} 
\dfrac{\partial \phi}{\partial y} (X_\tau,Y_\tau, \tau) \d W_\tau.
\]
Therefore $M_t^{\lambda,\phi}$ is a martingale with respect to $\mathcal{F}_t$. Similarly for $\phi'$,
\[
M_t^{\mu,\phi'}  = \phi'(x,y,0) + \int_0^t  e^{-\mu \tau} 
\dfrac{\partial \phi'}{\partial y} (X_\tau,Y_\tau, \tau) \d W_\tau.
\]
Hence, \eqref{ItoIsometry} is obtained using Ito's isometry. 
\end{proof}
\subsection{Proof for the backward-in-time parabolic problems.}
\begin{proof}
Applying \eqref{martingale1}-\eqref{martingale2} of Proposition \ref{prop1} to $\phi = u$, the solution of \eqref{cu}, we get
\[
\mathbb{E} \left ( \Gamma_{T}^\lambda(f,g) \right ) = u(x,y,0).
\]
Applying \eqref{ItoIsometry} of Proposition \ref{prop1} to $\phi = u$ and $\phi = v$, we get
\[
\mathbb{E} \left ( \Gamma_{T}^\lambda(f,g) \Gamma_{T+h}^\mu(\varphi,\psi) \right )
= uv(x,y,0) + 
\mathbb{E} \left ( \int_0^{T} e^{-(\lambda + \mu) \tau} \frac{\partial u}{\partial y} \frac{\partial v}{\partial y} (X_\tau,Y_\tau,\tau) \d \tau \right ).
\]
Finally, applying \eqref{martingale1}-\eqref{martingale2} of Proposition \ref{prop1} to $\phi = w$, the solution of \eqref{cw}, we get
\[
w(x,y,0) = \mathbb{E} \left ( \int_0^{T} e^{-(\lambda + \mu) \tau} \frac{\partial u}{\partial y} \frac{\partial v}{\partial y} (X_\tau,Y_\tau,\tau) \d \tau \right ).
\]
Therefore, 
\[
\mathbb{E} \left ( \Gamma_{T}^\lambda(f,g) \Gamma_{T+h}^\mu(\varphi,\psi) \right )  
= \left (uv + w \right )(x,y,0).
\]
\end{proof}
\subsection{Proof for the degenerate elliptic problems}
\begin{proof}
Under Assumption \eqref{A1}, using Proposition \ref{prop1},
\begin{align*}
& u_\lambda(x,y) = \mathbb{E} \int_0^\infty \exp(-\lambda t) g(X_t,Y_t) \d t, \quad v_\mu(x,y) = \mathbb{E} \int_0^\infty \exp(-\mu t) \psi(X_t,Y_t) \d t\\
& w_{\lambda+\mu}(x,y) =  \mathbb{E} \int_0^\infty \exp(-(\lambda+\mu) t) \frac{\partial u_\lambda}{\partial y}(X_t,Y_t) \frac{\partial u_\mu}{\partial y}(X_t,Y_t) \d t\\
\end{align*}
and under Assumption \eqref{A2}, direct calculations yields
\begin{align*}
& \mathbb{E} \int_0^\infty \exp(-(\lambda+\mu) t) \left ( \frac{\partial u_\lambda}{\partial y} \frac{\partial u_\mu}{\partial y} \right )(X_t,Y_t) \d t
+ \mathbb{E} \int_0^\infty \exp(-\lambda t) g(X_t,Y_t) \d t  \mathbb{E} \int_0^\infty \exp(-\mu t) \psi (X_t,Y_t) \d t\\
& =  \mathbb{E}  \int_0^\infty \int_0^\infty \exp(-\lambda s_1) \exp(-\mu s_2) g(X_{s_1},Y_{s_1}) \psi (X_{s_2},Y_{s_2}) \d s_1 \d s_2.
\end{align*}
\end{proof}

\section*{Acknowledgement}
\noindent \lmxxx{We thank the editor and the two anonymous reviewers for their helpful comments leading to the improvement of the present manuscript.} 
LM expresses his sincere gratitude to the Courant Institute
for being supported as Courant Instructor in 2014 and 2015, when this
work was initiated.  LM is also supported by a faculty discretionary
fund from NYU Shanghai and the National Natural Science Foundation of
China, Research Fund for International Young Scientists under the
project \#1161101053 entitled ``Computational methods for non-smooth
dynamical systems excited by random forces'' and the Young Scientist
Program under the project \#11601335 entitled ``Stochastic Control
Method in Probabilistic Engineering Mechanics''. LM also thanks the
Department of Mathematics of the City University of Hong Kong for the
hospitality. JW acknowledges support from the SAR Hong Kong grant 
[CityU 11306115] ``Dynamics of Noise-Driven Inelastic Particle Systems''.


\end{document}